\documentclass[11pt]{article}
\usepackage[dvipsnames]{xcolor}
\usepackage[utf8]{inputenc}
\usepackage{tikz}
\usetikzlibrary{shapes,arrows,positioning}
\usetikzlibrary{automata,arrows,positioning,calc}
\usepackage{amsmath}
\usepackage{amssymb}
\usepackage{amsthm}

\usepackage{dutchcal}
\usepackage{hyperref}

\usepackage{bbm}
\usepackage{blkarray}
\usepackage{caption}
\usepackage{subcaption}
\usepackage{float}
\usepackage{stmaryrd}
\usepackage[normalem]{ulem}
\usepackage{enumitem}

\usepackage{fancyhdr}
\usepackage{amsopn}

\usepackage{booktabs}

\definecolor{burgundy}{rgb}{0.5,0.0, 0.13}

\hypersetup{colorlinks,
	citecolor=burgundy,
	citebordercolor=burgundy,
	linkcolor=burgundy,
	urlcolor=MidnightBlue
}

\usepackage{natbib}
\def\newblock{\ }%
\bibpunct[, ]{[}{]}{,}{n}{}{,}%

\providecommand{\keywords}[1]
{
	\small	
	\textbf{{Keywords.}} #1
}

\numberwithin{equation}{section}
\usepackage{todonotes}

\newtheorem{theorem}{Theorem}[section]
\newtheorem{remark}[theorem]{Remark}

\newtheorem{definition}[theorem]{Definition}

\DeclareMathOperator*{\argmax}{arg\,max}
\DeclareMathOperator*{\argmin}{arg\,min}

\usepackage{algorithm,algorithmicx}
\usepackage[noend]{algpseudocode}

\usepackage{mathtools}
\mathtoolsset{showonlyrefs=true}

\usepackage{empheq}

\def\balpha{\boldsymbol{\alpha}}

\def\bmu{\boldsymbol{\mu}}

\def\nprime{{n^\prime}}

\newcommand\bX{\boldsymbol{X}}
\newcommand\bY{\boldsymbol{Y}}
\newcommand\bZ{\boldsymbol{Z}}

\newcommand\bW{\boldsymbol{W}}
\newcommand\bu{\boldsymbol{u}}
\newcommand\bm{\boldsymbol{m}}
\newcommand\bS{\boldsymbol{S}}

\newcommand\EE{\mathbb E}

\newcommand\RR{\mathbb R}

\newcommand\bbA{\mathbb A}

\newcommand{\indic}{\mathbf{1}}

\title{Machine Learning Methods for Large Population Games \\with Applications in Operations Research}

\author{G\"ok\c ce Dayan{\i}kl{\i}\footnote{Department of Statistics, 
	University of Illinois at Urbana-Champaign, Champaign, IL 61820, USA
		\href{mailto:gokced@illinois.edu}{gokced@illinois.edu}.}
	\and Mathieu Lauri\`ere\footnote{Shanghai Frontiers Science Center of Artificial Intelligence and Deep Learning; NYU-ECNU Institute of Mathematical Sciences, NYU Shanghai, 567 West Yangsi Road, Shanghai, 200126, People’s Republic of China, \href{mailto:mathieu.lauriere@nyu.edu}{mathieu.lauriere@nyu.edu}.}
}
\date{}

\begin{document}

\maketitle

\begin{abstract}
In this tutorial, we provide an introduction to machine learning methods for finding Nash equilibria in games with large number of agents. These types of problems are important for the operations research community because of their applicability to real life situations such as control of epidemics, optimal decisions in financial markets, electricity grid management, or traffic control for self-driving cars. We start the tutorial by introducing stochastic optimal control problems for a single agent, in discrete time and in continuous time. Then, we present the framework of dynamic games with finite number of agents. To tackle games with a very large number of agents, we discuss the paradigm of mean field games, which provides an efficient way to compute approximate Nash equilibria. Based on this approach, we discuss machine learning algorithms for such problems. First in the context of discrete time games, we introduce fixed point based methods and related methods based on reinforcement learning. Second, we discuss machine learning methods that are specific to continuous time problems, by building on optimality conditions phrased in terms of stochastic or partial differential equations. Several examples and numerical illustrations of problems arising in operations research are provided along the way.
\end{abstract}

\vskip3mm

\maketitle

\tableofcontents

\section{Introduction}

\noindent {\bf Motivations. } Optimization is pervasive in operations research (OR), and cover some classical examples such as supply and demand management in fast moving consumer products sector, healthcare operations management, or shortest path design for package distribution; see e.g.~\citet{rardin1998optimization}, \citet{frazier2018bayesian}, \citet{gill2019practical}. Optimal control theory provides a framework and a set of tools which extends standard optimization to dynamical problems. Stochastic optimal control further extends the theory of deterministic control to situations where disturbances may affect the dynamics. Applications range from adjusting the temperature of a room to controlling the trajectory of an aircraft, or deciding an optimal portfolio to maximize financial return. Stochastic optimal control theory has been developed extensively over the twentieth century and various numerical methods have been introduced, such as dynamic programming which was introduced by~\citet{bellman1958dynamic}, which exploits the time structure in a crucial way. More recently, a class of discrete time optimal control problems called Markov decision processes (MDPs) have been the focus of reinforcement learning (RL) methods, borrowing ideas from dynamic programming or direct policy optimization; see e.g.~\citet{sutton2018reinforcement,bertsekas2019reinforcement} for monographs on this topic. 
However, stochastic optimal control is limited to optimizing the performance of a single dynamical system but many real-world scenarios involve multiple strategic agents simultaneously controlling systems that interact through their dynamics or their objective functions. Examples can be found in traffic routing (drivers interact through the creation of traffic jams), finance (investors interact through a collective impact on the price), electricity management (producers and consumers need to reach an equilibrium between offer and demand of electricity), or epidemic management (citizens interact through the spread of the infection), to cite just a few. When the agents are non-cooperative, each controller must try to anticipate other agents' decisions in order to optimize their own objective function. Game theory has been introduced to study such situations, first in the static setting by~\citet{nash1951non,von1928theory} and then in dynamic setting by~\citet{fudenberg1991game,bacsar1998dynamic}. Dynamic game theory is significantly more complex than optimal control theory but has a huge potential in the OR applications. One of the bottlenecks for real-world applications has thus far been the development of efficient numerical methods for large scale games, but new approaches have been introduced in the recent years. We refer e.g. to~\citet{hu2023recent} for a survey of machine learning methods for stochastic optimal control and games mostly in continuous time, and to~\citet{lauriere2022learning} for a survey of learning methods for discrete time mean field games. Compared with these references, the present tutorial will cover discrete and continuous time models, with a focus on potential applications in OR.

\vskip 6pt
\noindent {\bf Scope of the tutorial. } The goal of this tutorial is to present recent methods for games with many players. In this direction, we will more specifically focus on two aspects of scalability: in terms of number of agents and in terms of model complexity. First, in order to model situations with a large number of interacting agents, we will introduce the framework of mean field games, see e.g.~\citet{huang2006large,lasry2006jeux,lasry2006jeux2}, which borrows ideas from statistical physics to provide a tractable approximation of games with very large populations. MFG theory has attracted a growing interest since its introduction, and various models have been introduced for applications that could be relevant in OR such as applications in epidemic control, see e.g.~\citet{laguzet2015individual,epidemics_SMFG,arefin2020mean,elie2020contact,aurell2022finite}, in energy demand and climate change policy decisions, see e.g.~\citet{aid2020entry,djehiche2020price,carmona2022mean,dayanikli2023multi,alasseur2020extended,bagagiolo2014mean}, in traffic control, see e.g.~\citet{chevalier2015micro,huang2021dynamic,festa2018mean}, in cybersecurity, see e.g.~\citet{kolokoltsov2016mean,miao2017cyber,kolokoltsov2018corruption}, in advertisement decisions, see e.g.~\citet{carmona2021mean,salhab2022dynamic}, or in systemic risk in the financial markets, see e.g.~\citet{carmona2015mean}.  
Second, to be able to solve complex models (e.g., problems in which the state dimension is high), we will move beyond the classical numerical approaches and will focus on machine learning based methods. 
Although more traditional methods benefit from the rich background of numerical analysis, machine learning methods (in particular deep learning and reinforcement learning methods) have met impressive empirical successes, which make them very appealing for applications in game theory and OR.

\vskip 6pt
\noindent {\bf Structure of the tutorial. } The rest of the paper is organized as follows. In Section~\ref{sec:finite-mf-games}, we introduce basic notations and definitions for optimal control and dynamic games in stochastic environments, both in discrete and continuous time. We also present the framework of mean field games, which provide a way to approximate games with very large number of players. To make these concepts more concrete and to show the flexibility of the frameworks, Section~\ref{sec:examples-extensions} presents several examples of models and discusses a few extensions which make the models more realistic. 
We then present classical families of machine learning methods for discrete time games in Section~\ref{sec:methods-discrete-time}. Machine learning methods for continuous time games are presented in Section~\ref{sec:methods-continuous-time} by exploiting the connection with stochastic differential equations and partial differential equations. Finally, we conclude with a summary and some future directions in Section~\ref{sec:conclusion}.

\section{Finite player games and mean field games}
\label{sec:finite-mf-games}

In this section, we first briefly present the framework of optimal control for a single agent\footnote{Throughout this tutorial, we use terms \textit{agent} and \textit{player} interchangeably.}, in discrete time and continuous time. We then present the framework of dynamical games, again in discrete and continuous times. From there, we present the paradigm of mean field games.

\subsection{Background on single-agent control}
\label{subsec:singleagent_background}

We will start by introducing the stochastic optimal control problems both in discrete and continuous times which are similar but differ in notation and terminology. In these problems, there is only one agent who aims to optimize their objectives by choosing some controls (i.e. actions).
\vskip2mm
\noindent \textbf{Discrete time. } Because of their wide applicability in the reinforcement learning, we will mainly focus on the \textit{Markov Decision Processes} (MDPs) for the discrete time setup. An MDP is defined by a tuple $(S, A, N_T, \mu_0, P, r)$ where $S$ is the state space, $A$ is the action space, $N_T>0$ is the total number of time steps, $\mu_0 \in \mathcal{P}(S)$ is an initial state distribution, $P_n(s, a) \in \mathcal{P}(S)$ is the probability distribution of the next state when using action $a$ in state $s$ at time $n$, and $r_n(s, a)$ is the immediate reward acquired when the agent used action $a$ in state $s$ at time $n$. To alleviate notation, we will write $\mathcal{T} = \{0,1,\dots,N_T\}$ for the set of times. Sometimes the next state is also included in the reward function's inputs. The actions and states can be both continuous or finite. We mostly focus on the finite action and state space setup since in the presentation of the continuous time setup we will focus on the continuous action and state spaces. The aim of the agent is to find an optimal \textit{policy}. The policy is a function $\pi: \mathcal{T} \times S \rightarrow \mathcal{P}(A)$ that gives the probability distribution on the actions that can be taken while being at a specific state $s$ at a given time. If this distribution is a Dirac distribution, then the policy is called a pure policy. Then, the objective of the agent is to choose a policy $\pi$ that maximizes the cumulative expected reward:
\begin{equation*}
\begin{aligned}
    J(\pi) = \ & \EE\left[\sum_{n=0}^{N_T} r_n(s_n, a_n) \right]\\
    & a_n \sim \pi_n(s_n),\ s_{n+1}\sim P_n(s_n, a_n),\ n \ge 0,\ s_0 \sim \mu_0.
\end{aligned}
\end{equation*}
In general, the last reward does not depend on the action, but we keep the same notation for all time steps to alleviate the presentation. 
The problem is introduced in the finite time horizon setting for the sake of consistency with the following sections on optimal control and mean field games (which are most often studied in finite horizon). However, it can be also stated in the infinite time horizon setting, in which case one usually considers the discounted reward and looks for stationary policies because this is sufficient to achieve the optimal value. We refer e.g. to the book of~\citet[Section 9.5]{bertsekas1996stochastic} for the rigorous proof of existence of an optimal policy that is stationary, under suitable assumptions. The reader interested in RL is referred to the book of~\citet{sutton2018reinforcement} for more background on MDPs from the machine learning perspective. 
One of the most common approaches to solve MDPs is to use dynamic programming, which relies on the introduction of a value function. We discuss in more details this notion in Section~\ref{sec:BR-by-DP} below.

\vskip2mm

\noindent \textbf{Continuous time. } 
The agent chooses an $A$-valued\footnote{In order to show the existence and uniquenesss results the control set $A\subseteq\RR^k$ is generally assumed to be closed, convex and bounded.} square-integrable control $\balpha := (\alpha_t)_{t \in [0,T]} \in \bbA$ to minimize their expected cost over a time horizon $[0,T]$ where $T>0$. In some applications, $T$ can also be taken as infinity, which is called infinite horizon stochastic optimal control problem. Controls can be of several form, such as open-loop, closed-loop, or Markovian. For the sake of brevity, we omit technical details here and refer interested reader to e.g.~\citet[Section 2.2]{borkar2005controlled} for more details. The agent's state $\bX:= (X_t)_t$ where $X_t \in \RR^d$ will evolve according to a stochastic differential equation. The dimension of the state and control can be integers $d$ and $k$, respectively; however, for the sake of simplification in notation, we will take $d=k=1$. In this section, we will focus on continuous state; however, in general, state can also be finite (i.e. discrete). The stochastic optimal control problem can be written mathematically as follows: minimize over $\balpha = (\alpha_t)_t$ the total expected cost
\begin{equation}
\label{eq:continuous-time-optimal-control-pb}
\begin{aligned}
    J(\balpha) = \ & \mathbb{E}\Big[\int_0^T {f(t, X_t, \alpha_t)} dt +g(X_T)\Big]\\
    &  dX_t = b(t, X_t, \alpha_t)dt + \sigma(t, X_t) dW_t, \quad X_0 =\zeta \sim \mu_0, 
\end{aligned}
\end{equation}
where $\bW = (W_t)_t$ is one-dimensional Brownian motion to represent randomness. Here, $f:[0,T] \times \RR \times A \rightarrow \RR$ is called running cost, $g:\RR \rightarrow \RR$ is called terminal cost, $b:[0,T] \times \RR \times A \rightarrow \RR$ is called drift, and $\sigma:[0,T] \times \RR \rightarrow \RR$ is called volatility. By following convention in the engineering community, we wrote the model as a minimization problem; however, it can also be written as a maximization problem. 
There are different ways to solve this problem by characterizing the solution with some differential equations. One approach is to use Hamilton-Jacobi-Bellman (HJB) equation to characterize the value function, another one is to use a backward stochastic differential equation by using Pontryagin stochastic maximum principle (see e.g.~\citet{yong2012stochastic}). We can introduce an application example as follows: Assume there is an employee who decides on how much effort at time $t$ they are going to put in the project where the value of the project at time $t$ can be thought as the state of the agent at that time. The value of the project will depend on how much effort the employee puts in and some constant exogenous randomness. Objectives of the employee can be to minimize their effort in the project while maximizing their time dependent and terminal utility from the project value. Then this problem can be mathematically modeled by choosing $b(t, X_t, \alpha_t) = \alpha_t,\sigma(t, X_t)=\sigma,f(t, X_t, \alpha_t)=\alpha_t^2-U(X_t)$, and $g(X_T)=-U(X_T)$ where $U:\RR\rightarrow\RR$ is a utility function and $\sigma>0$ is a constant. Realize that maximizing the utility corresponds to minimizing negative utility which in turn explains the minus signs in front of the utility functions.

\begin{remark}
    The volatility $\sigma$ can be a function of the control too, i.e., of the form $\sigma(t, X_t, \alpha_t)$. However, this requires using the \textit{full} Hamiltonian and is technically more challenging, both from the theoretical and the numerical viewpoints. Our presentation of the stochastic control problem is motivated by the mean field game framework we will discuss in the sequel. We choose to follow immediately the setting used in the vast majority of the MFG literature, which focuses on \textit{uncontrolled} volatility (i.e., the volatility is not a function of control $\alpha_t$). 
\end{remark}

\subsection{Nash equilibria in finite player games}

\noindent {\bf Discrete time. }
A discrete time $N$-player game is defined by a tuple $(\underline{S}, \underline{A}, N_T, \underline{\mu}^0, \underline{P}, \underline{r}) = ((S^1,\dots,S^N),$\\$ (A^1,\dots,A^N), N_T, (\mu_0^1,\dots,\mu_0^N), (P^1,\dots,P^N), (r^1,\dots,r^N))$ where $N$ is the number of players and the other symbols roughly have the same interpretation as before for each player, but now $P^i$ and $r^i$ are going to involve interactions between the players. $S^i$ and $A^i$ are the state and action spaces for player $i$, $\mu_0^i$ is the initial state distribution for player $i$, $P^i: \mathcal{T} \times \underline{S} \times A^i \to \mathcal{P}(S^i)$, $r^i: \mathcal{T} \times \underline{S} \times A^i \to \mathbb{R}$, and $P^i_n(\underline{s},a^i)$ and $r^i_n(\underline{s},a^i)$ denote respectively the probability distribution of $s_{n+1}^i$ and the reward for player $i$, when the states of all players are $\underline{s} = (s^1,\dots,s^N)$ and the action of player $i$ is $a^i$ at time $n$. A policy for player $i$ is a function $\pi^i: \mathcal{T} \times \underline{S} \rightarrow \mathcal{P}(A^i)$ that gives the probability distribution on the actions that can be taken while the players are at states $\underline{s}$ at a given time. Then, given the policies $\underline{\pi}^{-i} = (\pi^1,\dots,\pi^{i-1}, \pi^{i+1},\dots,\pi^N)$ used by other players, the objective of the agent $i$ is to choose a policy $\pi^i$ that maximizes the cumulative expected reward:
\begin{equation}
\label{eq:discrete-time-finite-player-reward}
\begin{aligned}
    J^i(\pi^i; \underline{\pi}^{-i})  =& \EE\left[\sum_{n=0}^{N_T} r_n^i(\underline{s}_n, a^i_n) \right],\\
    & a^j_n \sim \pi^j_n(\underline{s}_n),\ s^j_{n+1}\sim P^j_n(\underline{s}_n, a^j_n),\ n \ge 0, s^j_0 \sim \mu^j_0, j \in [N].
\end{aligned}
\end{equation}
Here again, the problem could be studied in the infinite-horizon discounted setting. Since the policy of player $i$ depends on the policy of the other players, the problem is, thus far, not fully well-defined. The most common notion of solution is Nash equilibrium, in which each player has no incentive to deviate unilaterally. More precisely:

\begin{definition}
\label{def:discrete-time-Nash}
    A Nash equilibrium is a policy profile $\underline{\hat\pi} = (\hat\pi^1,\dots,\hat\pi^N)$ such that, for every $i \in [N]$, $\hat\pi^i$ is a maximizer of $J^i(\cdot; \underline{\hat\pi}^{-i})$. 
\end{definition}

In the sequel, to be able to pass to the limit when the number $N$ of players goes to infinity, we will focus on the case where the players are indistinguishable and have symmetric interactions. In other words, we will assume $S^i=S, A^i=A$, $\mu_0^i = \mu_0$, $P^i = P$ and $r^i=r$ for all $i$, and furthermore, for every $n$ and $a^i$, $P_n(\cdot,a^i)$ and $r_n(\cdot,a^i)$ are symmetric functions of the input $\underline{s} = (s^1,\dots,s^N)$, i.e., they are invariant under permutation of the $N$ coordinates of $\underline{s}$. This structure is crucial to be able to pass to the mean field limit, as we will discuss in the next section.

\noindent {\bf Continuous time. } 
In continuous time, each agent solves an optimal control problem analogous to~\eqref{eq:continuous-time-optimal-control-pb} except that, in the game setting, the cost and the dynamics potentially depend on the other players states or actions. For simplicity of presentation, we will only consider interactions through states.

Let us consider $N>1$ agents. Agent $i$ chooses an $A^i$-valued square-integrable control $\balpha^i := (\alpha_t^i)_t \in \bbA^i$ to minimize their expected cost over a time horizon $[0,T]$ where $T>0$. The agent's state $\bX^i := (X^i_t)_t$ where $X^i_t \in \RR^d$ will evolve according to a stochastic differential equation (SDE). As before, for the simplicity in presentation, we will take the state and action dimensions equal to 1 and we will focus on continuous state and control spaces. Similar to the continuous time stochastic control models, the state and the control spaces can be also finite. Given the controls of the other players $\underline\balpha^{-i} = (\balpha^1, \dots, \balpha^{i-1}, \balpha^{i+1}, \balpha^N)$, the goal for player $i$ is to solve the stochastic optimal control problem: Minimize over $\balpha^i$ 
\begin{equation}
\label{eq:continuous-time-optimal-control-pb-player-i}
\begin{aligned}
    J(\balpha^i; \underline\balpha^{-i})= \, & \mathbb{E}\Big[\int_0^T {f^i(t, \underline{X}_t, \alpha^i_t)} dt +g^i(\underline{X}_T)\Big]\\
    &  dX^i_t = b^i(t, \underline{X}_t, \alpha^i_t)dt + \sigma^i(t, \underline{X}_t) dW^i_t, \quad X^i_0 =\zeta^i \sim \mu^i_0, 
\end{aligned}
\end{equation}
where $\bW^i = (W^i_t)_t$ is one-dimensional Brownian motion to represent randomness specific to player $i$. The random variables $\zeta^j$ and the Brownian motions $W^j$ are all independent. Here, the running cost $f^i:[0,T] \times \RR^N \times A \rightarrow \RR$, the terminal cost $g^i:\RR^N \rightarrow \RR$, the drift $b^i:[0,T] \times \RR^N \times A \rightarrow \RR$, and the volatility $\sigma^i:[0,T] \times \RR^N \rightarrow \RR$ depend on the positions of all the players. They could also depend on the actions of all the players but we omit this to alleviate the notations. 

In this context, a Nash equilibrium is defined as follows.

\begin{definition}
    A Nash equilibrium is a control profile $\underline{\hat\balpha} = (\hat\balpha^1,\dots,\hat\balpha^N)$ such that, for every $i \in [N]$, $\hat\balpha^i$ is a minimizer of $J^i(\cdot; \underline{\hat\balpha}^{-i})$. 
\end{definition}

As already mentioned in the discrete time case, to pass to the mean field limit when $N$ goes to infinity, we will focus on the case where the players are indistinguishable and have symmetric interactions. In the context of the above formulation in continuous time, we will assume $\mu_0^i = \mu_0$, $b^i = b$, $\sigma^i = \sigma$, $f^i = f$, and $g^i = g$ for all $i$, and furthermore, for every $t$ and $a^i$, $f(t,\cdot,a^i)$, $g(\cdot)$, $b(t,\cdot,a^i)$ and $\sigma(\cdot)$ are symmetric functions of the input $\underline{x} = (x^1,\dots,x^n)$, i.e., they are invariant under permutation of the $N$ coordinates of $\underline{x}$.

\subsection{From large populations to mean field games}

\begin{figure}[H]
\centering
    \includegraphics[width=0.6\textwidth]{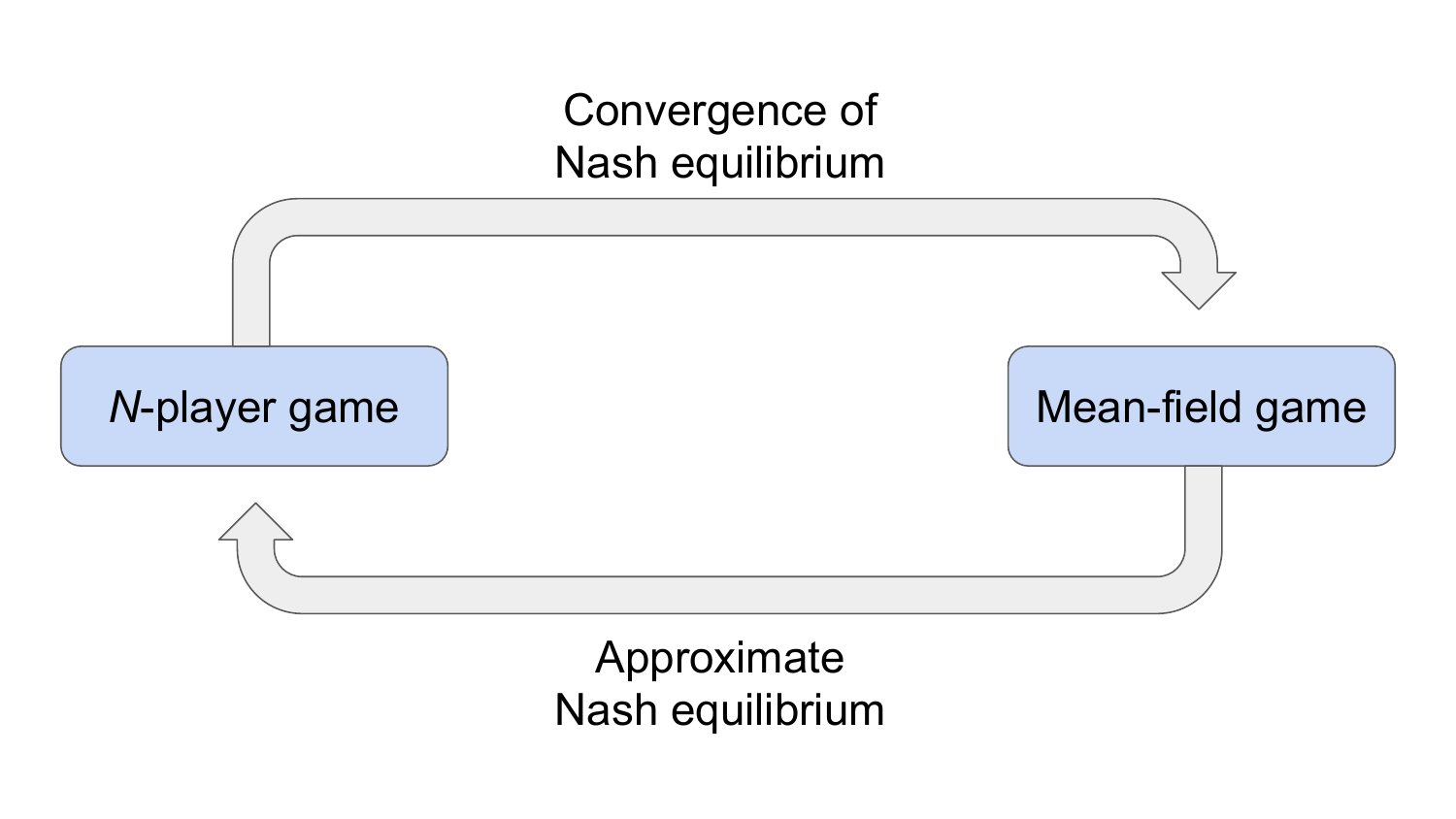}
\caption{Conceptual connections between $N$-player game and mean field game. Under suitable assumptions, letting $N$ go to infinity (from left to right), the Nash equilibria in an $N$-player game converge to the mean field Nash equilibrium. Conversely (from right to left), using a mean field Nash equilibrium in a finite player game leads to an approximate Nash equilibrium whose quality improves with $N$.} 
\label{fig:diagram-N-player-MF}
\end{figure}

As the number of agents gets higher, we need to consider an exponentially increasing number of interactions which will create difficulty in tractability. In order to prevent this, mean field games (MFGs) have been recently proposed to approximate the games with large number of agents. In this section, we will introduce MFGs both in discrete and continuous time setups and later on, solving these models with machine learning will be the focus of this tutorial. In MFGs, all the agents are assumed to be identical and interacting symmetrically through the population distribution (of state, action/control, or joint state-control) and the number of agents is assumed to go to infinity. Thanks to these assumptions, we can focus on a \textit{representative} agent who is infinitesimal such that when she changes her control, the population distribution is not affected. This reduces the problem of $N$ interacting agents to a problem of a representative agent and her interactions with the population distribution.
 At a high level, the connection between finite player games and mean field games are summarized in Figure~\ref{fig:diagram-N-player-MF}. We discuss this diagram in more details at the end of this section.

\vskip2mm
\noindent \textbf{Discrete time. } 
In an MFG, we want to understand the evolution of the mean field (which in our case will be the state distribution) and the decisions made by one representative player. This is sufficient to understand the behavior of the whole population since all the players are assumed to be indistinguishable.

A discrete time MFG is defined by a tuple $(S, A, N_T, \mu^0, P, r)$ where the symbols have the same interpretation as before, but now $P$ and $r$ are going to involve interactions with the mean field. $S$ and $A$ are the state and action spaces for each player, $\mu_0$ is the initial state distribution for each player, $P: \mathcal{T} \times S \times A \times \mathcal{P}(S) \to \mathcal{P}(S)$, $r: \mathcal{T} \times S \times A \times \mathcal{P}(S) \to \mathbb{R}$, and $P_n(s,a,\mu)$ and $r_n(s,a,\mu)$ denote respectively representative player's next state probability distribution and reward, if at time $n$ this player is at state $s$, uses action $a$ and the rest of the population has state distribution $\mu$. A policy for a player is a function $\pi: \mathcal{T} \times \underline{S} \rightarrow \mathcal{P}(A)$ that gives the probability distribution on the actions that can be taken at a given time. In this formulation, the policy is a function of the representative player's own state and not a function of the mean field distribution. The reason for this comes from the fact that when the number of players goes to infinity, the mean field will be a deterministic function of time through the law of large numbers. Therefore, its effect can be taken into account directly through the time dependence of the policy.

 It will be convenient to use the following notation: given a policy $\pi$, we denote by $\mu^\pi$ the mean field sequence induced when the whole population uses policy $\pi$, i.e., 
 \begin{equation}
 \label{eq:discrete-time-mean-field-evol}
    \mu^\pi_0 = \mu^0, \quad \mu^\pi_{n+1}(s) = \sum_{s^\prime,a} \mu^\pi_n(s^\prime) \pi_n(a|s^\prime) P_n(s^\prime,a,\mu^\pi_n),  n=0,\dots,N_T-1.
\end{equation}
 
 Here, we focus on interactions through the state distribution only, although more general interactions (e.g. through the action distribution) could be considered. As a consequence, a player's reward can be defined as a function of the mean field sequence instead of a function of the other player's policies. This is a major simplification compared with finite-player games (see~\eqref{eq:discrete-time-finite-player-reward}). Another key simplification is that each player has no influence on the mean field at all, so unilateral deviations do not perturb the mean field. Optimality conditions are thus easier to phrase than in finite-player games, as we shall see in the sequel.   
 Given a mean field $\mu = (\mu_n)_{n \in \mathcal{T}}$ describing the behavior of the rest of the population, the objective of a representative player is to choose a policy $\pi$ that maximizes the cumulative expected reward:
\begin{equation}
\label{eq:discrete-time-mfg-reward}
\begin{aligned}
    J(\pi; \mu)  = \, & \EE\left[\sum_{n=0}^{N_T} r_n(s_n, a_n, \mu_n) \right],\\
    & a_n \sim \pi_n(s_n), s_{n+1}\sim P_n(s_n, a_n, \mu_n), n \ge 0, s_0 \sim \mu_0.
\end{aligned}
\end{equation}
Based on this notion, we can define the concept of Nash equilibrium. 
The counterpart to Definition~\ref{def:discrete-time-Nash} is:
\begin{definition}
\label{def:MFGNE-discrete}
    A (discrete time mean field) Nash equilibrium is a policy $\hat\pi$ and a mean field $\hat\mu$ such that the following two conditions are satisfied:
    \begin{itemize}
        \item $\hat\pi$ is a maximizer of $J(\cdot; \hat\mu)$ as defined in~\eqref{eq:discrete-time-mfg-reward},
        \item $\hat\mu$ is the mean field induced by policy $\hat\pi$, i.e., $\hat\mu = \mu^{\hat\pi}$ using the notation introduced in~\eqref{eq:discrete-time-mean-field-evol}.
    \end{itemize}
\end{definition}
This definition can be interpreted as a fixed point problem. Let us denote by $BR$ the function that gives the set of best responses given a mean field, and by $MF$ the function that gives the mean field induced by a policy. Then a mean field Nash equilibrium $(\hat\pi,\hat\mu)$ satisfies $\hat\pi \in BR \circ MF(\hat\pi)$ and $\hat\mu \in MF \circ BR(\hat\mu)$ where $\circ$ denotes the composition of functions. So solving an MFG amounts to finding a fixed point of $BR \circ MF$ or $MF \circ BR$. It should be noted that, while uniqueness of the equilibrium mean field $\hat\mu$ can be guaranteed by monotonicity assumptions, it is in general more difficult to ensure uniqueness of the best response to $\hat\mu$. In the continuous space setting, uniqueness of the best response can be guaranteed under strict convexity assumptions. The existence of the discrete time mean field Nash equilibrium can be proved either: (a)  by showing that the of the fixed point mappings ($BR \circ MF$ or $MF \circ BR$) is a strict contraction and using Banach fixed point theorem, or (b) by proving some compactness property and using Schauder or Kakutani fixed point theorem for instance, see e.g.~\citet{gomes2010discrete}. The first method also gives uniqueness by typically requires quite strong assumptions (e.g., Lipschitz continuity of the dynamics and cost with small enough Lipschitz constants).

It can be shown, under suitable conditions, that a mean field Nash equilibrium policy $\hat\pi$ is an approximate Nash equilibrium for the corresponding finite player game, and the quality of the approximation improves as the number of player increases (see at the end of the section for more details). Furthermore, it is in general computationally more efficient to solve an MFG than to tackle directly the corresponding game with a finite but large number of players, see Figure~\ref{fig:computation_time_algo_n_vs_mean_field} below for an illustration in a routing game.

\vskip2mm
\noindent \textbf{Continuous time. } Similar to the previous sections, the \textit{representative} agent chooses an $A$-valued\footnote{As mentioned previously in Section~\ref{subsec:singleagent_background}, in order to show the existence and uniquenesss results the control set $A\subseteq\RR^k$ is generally assumed to be closed, convex and bounded.} square-integrable control $\balpha := (\alpha_t)_t \in \bbA$ to minimize their expected cost over a time horizon $[0,T]$ where $T>0$ (as before $T$ can be infinity). The representative agent's state $\bX:= (X_t)_t$ where $X_t \in \RR^d$ will evolve according to an SDE that depends on her own state, control, as well as the population distribution. For the simplicity in presentation, we will take the dimension of state and control equal to 1. As in the previous sections, we will present the model with continuous state space; however, in general, state can also be finite (i.e. discrete). For example, we can have an application related to mitigation of epidemics. In this case, the agents' states will be their health status such as susceptible, infected, recovered, etc. The objective of the representative agent using control $\balpha\in\bbA$ when interacting with other agents through the population state distribution, $\bmu:= (\mu_t)_t$ with $\mu_t\in \mathcal{P}(\RR)$ is to minimize
\begin{equation}
\label{eq:continuous-time-mfg-cost}
\begin{aligned}
    J(\balpha; \bmu) = \, &\mathbb{E}\Big[
    \int_0^T f(t, X_{t},\alpha_t, \mu_t)dt  + g(X_T, {\mu_T})\Big]
    \\
    & dX_t =  b(t, X_t, \alpha_t, \mu_t)dt + \sigma(t, X_t, \mu_t) dW_t, \qquad X_0 =\zeta \sim \mu_0,
\end{aligned}
\end{equation}
where $\bW := (W_t)_t$ is the Brownian Motion representing the idiosyncratic (i.e. individual and independent) noise the agents face. The MFG models where there is a common randomness source all the agents are affected in the same way are called \textit{MFGs with common noise} which are more technically involved since the players' states are not independent anymore and one needs to work with a conditional mean field. Following the terminology used in stochastic optimal control problems, we call function $f, g, b, \sigma$ running cost, terminal cost, drift and volatility, respectively. In the model above, the representative agent interacts with the population through the population state distribution; however, in general the interactions can be through the population control distribution or joint control and state distribution which is called \textit{extended} MFG models in the literature. The continuous time MFG models can be solved with the analytic approach or the probabilistic approach. In the former, the solution is characterized by forward-backward partial differential equations (FBPDE) where the backward component is the Hamilton-Jacobi-Bellman (HJB) equation and the forward component is the Kolmogorov-Fokker-Planck (KFP) equation. In the latter, the solution is characterized by forward-backward stochastic differential equations (FBSDE). We refer to the monographs of~\citet{bensoussan2013mean,carmona2018probabilistic} for more background.

In this context, we can define the notion of Nash equilibrium as follows. 
\begin{definition}
\label{def:MFGNE-continuous}
    A (continuous time mean field) Nash equilibrium is a control $\hat\balpha$ and a mean field $\hat\bmu$ such that the following two conditions are satisfied:
    \begin{itemize}
        \item $\hat\balpha$ is a minimizer of $J(\cdot; \hat\bmu)$ as defined in~\eqref{eq:continuous-time-mfg-cost},
        \item $\hat\bmu$ is the mean field induced by control $\hat\balpha$. 
    \end{itemize}
\end{definition}
As in the discrete time case, the problem can be phrased as a fixed point problem. This fixed point formulation provides a basis for existence results using fixed point theorems (see e.g.~\citet{huang2006large,carmona2018probabilistic}). Uniqueness of the equilibrium mean field can be proved under Lasry-Lions monotonicity condition, see e.g.~\citet{lasry2007mean}, while uniqueness of the best response is generally ensured by strict convexity assumptions that result in having a unique optimizer for the Hamiltonian.  

The two directions illustrated by Figure~\ref{fig:diagram-N-player-MF} have been proved rigorously in several settings, under suitable assumptions. For the justification that mean field Nash equilibrium provides a good approximate Nash equilibrium in finite player games, interested readers can refer for instance to~\citet{huang2006large,bensoussan2013mean} in the continuous setting, and to~\citet{saldi2019approximate,saldi2020approximate,cui2024learning} in the discrete setting. For the convergence of Nash equilibria in finite-player games towards mean field Nash equilibria, we refer e.g. to the works of~\citet{cecchin2019convergence,cardaliaguet2019master,lacker2020convergence}.

\section{Examples and extensions}
\label{sec:examples-extensions}
In this section, we provide several examples from the mean field games literature in both discrete time and continuous time. We conclude by mentioning several extensions of the framework.

\subsection{Discrete time models}

We first note that, in the same way that optimal control contains optimization problems as a special case, static (single-shot) games can naturally be embedded in the framework of finite horizon dynamic games by considering games with a single state and a single time step. We refer to e.g.~\citet{muller2022learningPSRO,muller2022learning,wang2023empirical} for variants of the classical rock-paper-scissors game and prisoner's dilemma game in the finite-horizon MFG framework, as well as learning algorithms for such games and refer to~\citet{carmona2021mean} for a static mean field game application. In this section, we present a few models which are truly dynamic.

\subsubsection{Crowd motion}
\label{sec:ex-crowd}
We start with an example in discrete time and discrete spaces which has a straightforward physical interpretation: crowd motion. Pedestrians move around in a spatial domain and interact through the density of the crowd, which may for instance create discomfort or slower movements. The main ingredients can be summarized as follows:
\begin{itemize}
    \item Agents: pedestrians.
    \item States: possible locations, for instance a grid world of the form $S = \{0,\dots,s_1\} \times \{0,\dots,s_2\}$.
    \item Actions: possible moves; for instance up, left, down, right or stay, which corresponds to $A = \{(1,0), (0,-1), (-1,0), (0,1), (0,0)\}$
    \item Transitions: $P_n(s,a,\mu) = $ the distribution of $s+a+\epsilon_n$, where $\epsilon_n$ is a random variable taking values in $A$ which is interpreted as a random perturbation of the agent's movements.
    \item Rewards: $r_n(s,a,\mu)$ can include spatial preferences (e.g., of the form $-|s-s_{*}|$ where $s_*$ is a target state), cost of moving (e.g., $-|a|$), and mean-field interactions such as:
    \begin{itemize}
        \item attraction towards the mean (e.g., $-|x - \bar\mu|$ where $\bar\mu$ is the mean of $\mu$);
        \item crowd aversion (e.g., $-\mu(x)$, which penalizes the fact of being in a crowded location) ;
        \item congestion (e.g., $-\mu(x)|a|$, which penalizes the fact of moving through a crowded region).
    \end{itemize}
    Notice that the last two types of interactions can also be non-local, e.g., $-(\rho \star \mu)(x)$ or $-(\rho \star\mu)(x)|a|$ where $\rho$ is a kernel and $\star$ denotes the (discrete) convolution.
    \item Obstacles can be included in the model by adding a set of forbidden states and, when an agent is supposed to move to a forbidden state, she just stays at her current state.
\end{itemize}
Numerical illustrations are provided in Section~\ref{sec:numerics-crowd} for a model of the above type, solved using fictitious play and reinforcement learning. Similar models have been considered for instance by~\citet{perrin2022generalization,zaman2023oracle,algumaei2023regularization}.

\subsubsection{Traffic routing}
\label{sec:example traffic}

We then consider another application with a clear physical interpretation: traffic routing. Cars evolve on a network of roads and try to reach a destination, but each car faces traffic jams created by the population of cars. The main elements of the model are:

\begin{itemize}
    \item Agents: cars (or drivers).
    \item States: one possible approach is to consider that the states are the locations along (a discretization of) the roads (edges of the network); however this could lead to a very large number of states; another approach is to consider that the state of a car is: the edge on the network and the waiting time before it can moves to the next edge.  
    \item Actions: in routing, the actions a car can take when it arrives at a cross-road is the next road (edge of the network); in a more sophisticated model, one could consider that the car can also choose its velocity (or acceleration) along the current road.   
    \item Transitions: they can model the movement of the car along the current road and then the transition to the next road; alternatively, they can model the decay of the waiting time until reaching the crossroad and then the transition to the next road; 
    \item Rewards: each wants to reach a target destination as soon as possible.
    \item If there are several origins and destinations, this can be represented using a multi-population MFG.
\end{itemize}
A model of this type has been considered for instance by~\citet{cabannes2022solving}. A more precise description of the model is given in Section~\ref{sec:numerics-routing} and numerical illustrations are provided.

\subsubsection{Cybersecurity}
Next, we present a model of cybersecurity introduced by~\citet{MR3575619}; see also the book of~\citet[Section 7.2.3]{MR3752669}. The players are computers which do not want to be infected, and there are finitely many states, corresponding to the level of infection and the level of protection. The original model is in continuous time but here we will present a discrete time version, to be consistent with our notations. More generally, continuous-time finite-state games can be approximated by discrete-time finite-state games by discretizing time and then replacing transition rate matrices by transition probabilities.

\begin{itemize}
    \item Agents: computers.
    \item States: there are 4 possible states, $S = \{DI, DS, UI, US\}$, which corresponds to the level of defense (Defended or Undefended) and the level of infection (Infected or Susceptible).
    \item Actions: each computed can decide that they want to change their level of protection, which is represented by the choice of a level $a$ in $A = \{0,1\}$; $a=0$ (resp. $a=1$ ) means that they want to keep (resp. switch to the other extreme) their current level of protection.
    \item Transitions: at a high level, the probability of being infected increases with the proportion of infected computers, the intensity of attack (which is a fixed parameter in this model) and it is lower when the computer is defended than when it is not. If one starts from the original continuous time model introduced by~\citet{MR3575619} with transition rate matrix denoted by $Q$, then the transition matrix $P$ over a time step of size $\Delta t$ is given by $P = \exp(\Delta t Q)$. The transition rate matrix in~\cite[Section 7.2.3]{MR3752669} is:
    $$
    	Q 
    	= \begin{pmatrix}
    	\dots 	& 		Q_{DS \rightarrow DI}	&	 \rho a 	&	0
    	\\
    	q_{rec}^D 	& 	\dots 		&	 0	&	\rho a
    	\\
    	\rho a 	& 	0 		&	 \dots	&	Q_{US \rightarrow UI}
    	\\
    	0	&	\rho a	&	q_{rec}^U	& \dots
    	\end{pmatrix}
    $$
    where 
    $
    	Q_{DS \rightarrow DI} = v_H q_{inf}^D + \beta_{DD} \mu(\{DI\})  + \beta_{UD} \mu(\{UI\}) ,
    $ and $Q_{US \rightarrow UI} = v_H q_{inf}^U+$ \\$ \beta_{UU} \mu(\{UI\}) + \beta_{DU} \mu(\{DI\}),
    $ 
    and all the instances of $\dots$ should be replaced by the negative of the sum of the entries of the row in which $\dots$ appears on the diagonal. The matrix depends on the action $a$ and the mean field $\mu$. The parameters are as follows: $\rho >0$ determines the rate at which a computer can switch its level of protection; $q_{rec}^D$ or $q_{rec}^U$ are the rates at which a computer can recover depending on whether it is defended or not; $v_H q_{inf}^D$ (resp. $v_H q_{inf}^U$) determines the rates at which a computer becomes infected by the attacker ($v_h$ represents the strength of the attack) if it is defended (resp. undefended); $\beta_{UU}\mu(\{UI\})$ (resp. $\beta_{UD}\mu(\{UI\})$) represents the rate of infection by undefended infected computers if the computer is undefended (resp. defended); $\beta_{DU}\mu(\{DI\})$ (resp. $\beta_{DD}\mu(\{DI\})$) represent the rates of infection by defended infected computers if it is undefended (resp. defended).
    
    \item Rewards: the rewards discourage the computers from being infected but also makes them pay a cost to be defended; it takes the form: $r(x,a,\mu) = -\left[ k_D \indic_{\{DI, DS\}}(x) + k_I \indic_{\{DI, UI\}}(x)\right]$.
\end{itemize}
Models of this type have been considered for instance by~\citet{lauriere2021numerical,vasal2023sequential}. 

\subsection{Continuous time models}

Some of the discrete time models can also be studied in continuous time. For example, crowd motion has been studied in the MFG literature by~\citet{lachapelle2011mean,achdou2015system,achdou2019mean}, while the cybersecurity model discussed above was introduced by~\citet{MR3575619} in continuous time as already mentioned, and then revisited e.g. by~\citet{kolokoltsov2018corruption,MR3752669}. We refer to these works for more details on these models. Below, we present several other models.

\subsubsection{Project value management} 
\label{cont_ex:employee}
We extend the example given in the introduction of continuous time stochastic optimal control problems in Section~\ref{subsec:singleagent_background}. In the extended model, we have infinitely many non-cooperative agents (i.e., employees) and the representative agent's problem is given as follows:
\begin{equation*}
\begin{aligned}
    \min_{(\alpha_t)_t}\ & \mathbb{E}\Big[\int_0^T \left(\frac{1}{2}\alpha_t^2 -  U(X_t)\right) dt - U(X_T)\Big]\\
    & \qquad dX_t = (\alpha_t + \bar{X}_t) dt + \sigma dW_t, \quad X_0= \zeta \sim \mu_0, 
\end{aligned}
\end{equation*}
where $U:\RR \rightarrow \RR$ is a concave utility function. Intuitively, the mathematical model implies that each player has some utility from their own project value (i.e., state) which is denoted by the $\RR$-valued process $\bX=(X_t)_{t\in[0,T]}$ and has some cost from putting too much effort (i.e., control) denoted by $\RR$-valued process $\balpha=(\alpha_t)_{t\in[0,T]}$. The initial project value of the representative agent follows an initial distribution $\mu_0$. We see that the dynamics of the value of the project is affected by the agent's own effort ($\alpha_t$) and also the average value of the projects in the population which is denoted by $\bar{X}_t=\int_{\RR} x  \mu_t(dx)$ at a given time. In this model, the representative player is interacting with the population through the average state (i.e., average project value) in the population. 

\subsubsection{Electricity production}
\label{cont_ex:electricity}
Motivated by the example given in~\citet*{carmona2022mean}, we present a game problem for large number of non-cooperative the electricity producers who are using nonrenewable energy resources in the production and who are paying a carbon tax depending on their pollution levels. The mean field interactions will come through the price of the electricity. In this setup, the representative electricity producer's model is given as follows:
\begin{equation*}
\begin{aligned}
    \min_{(\alpha_t)_t}\ & \mathbb{E}\Big[\int_0^T \left(c_1 (\alpha_t)^2 + p N_t +c_2 (Q_t-D_t)^2 - c_3(\rho_1 +\rho_2 (D_t-\bar{Q}_t)) Q_t \right) dt +\tau P_T\Big]\\
        &dN_t = \alpha_t dt, \qquad N_0\sim \mu_{n,0}, \qquad\qquad dE_t = \delta \alpha_t dt   + \sigma_1 dW^1_t, \qquad E_0\sim \mu_{e,0},\\
        &dP_t = E_t dt, \qquad P_0\sim \mu_{p,0}, \qquad \qquad dQ_t = \kappa \alpha_t dt   + \sigma_2 dW^2_t,  \qquad \ Q_0\sim \mu_{q,0}.
\end{aligned}
\end{equation*}
In this model, the representative electricity producer controls the change in nonrenewable energy source usage at time $t$ and it is denoted by $\alpha_t\in \RR$. The representative producer has a 4-dimensional state: nonrenewable energy usage at time $t$ ($N_t \in \RR_+$), instantaneous carbon emission at time $t$ ($E_t \in \RR_+$), cumulative pollution at time $t$ ($P_t \in \RR_+$), and instantaneous electricity production at time $t$ ($Q_t \in \RR_+$). In the dynamics, $\delta, \kappa, \sigma_1, \sigma_2$ are positive coefficients, where $\delta$ and $\kappa$ denote the carbon emission level and production efficiency per unit nonrenewable energy usage, respectively. $(W_t^1)_{t\in[0,T]}$ and $(W_t^2)_{t\in[0,T]}$ two independent Brownian motions that represent idiosyncratic noises. In the objective, the first term represents the cost of ramping up or down the production too quickly in the facility, the second one is the price paid for instantaneous nonrenewable energy usage where $p$ denotes the price per unit nonrenewable energy usage. The third term is the penalty for matching the exogenous average electricity demand which is denoted by $(D_t)_{t\in[0,T]}$, and the fourth term is the revenue gained from the selling the electricity. Here $\rho_1 +\rho_2 (D_t - \bar{Q}_t)$ denotes the price per unit of electricity. We can see that the price is increasing if the average demand is more than the average supply ($\bar Q_t=\EE[Q_t]$) and is decreasing vice versa. This term is introducing the mean field interactions. Finally, the fifth term represents the carbon tax the producer pays at the terminal time $T$.

\subsubsection{Trading with price impact} 
\label{sec:example trading}

Another application example of continuous time mean field games is modeling the traders. Here, for the sake of definiteness, we follow the model introduced by~\citet{MR3805247} and we assume that there is a single stock being traded. In this models the agents are traders, or more specifically brokers who are in charge of liquidating their portfolio in a given time interval. In this setup, the representative trader's model is given as follows:
\begin{equation*}
    \begin{aligned}
         \min_{(\alpha_t)_t} - &\EE \left[- \phi \int_0^T |Q_t|^2 dt +  X_T + Q_TS_T - A |Q_T|^2  \right]\\
         & dS_t = \gamma \bar{\mu}_t  dt + \sigma dW_t, \qquad\qquad\ S_0 \sim \mu_{s,0}\\
        &dQ_t = \alpha_t dt, \qquad\qquad\qquad\qquad\  Q_0 \sim \mu_{q,0}\\
        &dX_t = -\alpha_t(S_t+\kappa \alpha_t)dt, \qquad\quad X_0\sim\mu_{x,0}.
    \end{aligned}
\end{equation*}
In this model, the representative trader controls their rate of buying and selling the stock at time t which is denoted by $\alpha_t \in \RR$. The state of the representative trader has three components: the price of the underlying stock $S_t$, which is common to all the players, the individual portfolio $Q_t$ (i.e., number of shares of the stock), the individual wealth $X_t$ (for example in dollars). Therefore, the state space is $\RR^3$. Notice the portfolio and the wealth can potentially take negative values if the trader needs to buy shares and is allowed to borrow money. The representative trader's control directly affects their portfolio and wealth. The individual wealth process ($\bX = (X_t)_{t\in[0,T]})$) is affected by a \textit{temporary} price impact while the stock price process ($\bS = (S_t)_{t\in[0,T]})$) is influenced by a \textit{permanent} price impact. Their scales are respectively determined by positive coefficients $\kappa$ and $\gamma$. The permanent price impact involves $\bar{\mu}_t$, which is the average trading rate (average action), i.e., $\bar\mu_t = \EE[\alpha_t]$. We refer to~\citet{MR3805247} for more explanations. In the objective, the first term and the last term are the cost of holding a large portfolio (either at time $T$ or any time before that). Here, $A$ and $\phi$ are positive constants that weigh the importance of these penalizations and can be interpreted as parameterizing the risk preference of the trader. The second and the third terms represents the terminal payoff which depends on the wealth and the value of the portfolio at time $T$. Since the objective is \textit{minimizing}, a negative sign is added before the expectation.

Since the model involves interactions through the distribution (in this mode through its mean) of actions, it is referred to as an \textit{extended} MFG or an MFG of controls. We present numerical illustrations in Section~\ref{sec:numerics-trading}.

The model is similar to the model considered by~\citet{MR3500455}, which can be solved using reinforcement learning and real data as e.g. by~\citet{leal2020learning}. However, the main difference of the above model is the fact that the price impact is here endogeneous since it is the result of the actions of all players, which makes the model more realistic than when $\bar\mu_t$ is given as an exogeneous term. Other similar models of optimal execution have been proposed by~\citet{MR3325272} in the weak formulation, and revisited in the book of~\citet[Sections 1.3.2 and 4.7.1]{MR3752669} in the strong formulation.

\subsection{Extensions}

To simplify the presentation, we restrict our attention to \emph{homogeneous} populations with \emph{symmetric} interactions, and to the notion of \emph{Nash equilibrium}. However real-world situations are often much more complex. Several extensions of this framework have been studied. For instance, multi-population MFGs consider MFGs in which there are several populations who interact, and each population is approximated by a mean field (see e.g.~\citet{huang2006large,feleqi2013derivationmulti,cirant2015multi,bensoussan2018meanmulti} for an analytical approach and~\citet[Section 7.1.1]{carmona2018probabilistic} for a probabilistic formulation). Graphon games go one step further towards heterogeneous interactions by considering that the continuum of agents have interactions that are weighted by a network structure (see e.g.~\citet{parise2023graphon,carmona2022stochasticgraphonstatic} in the static case and e.g.~\citet{caines2019graphon,caines2021graphon,aurell2022stochastic} in the dynamic case; see~\citet{cui2021learning,fabian2023learning} for learning methods).
Major-minor MFGs include one (or several) significant (i.e., major) agent(s) to the model, where she (they) can directly affect the mean field (see e.g.~\citet{huang2010large, nguyen2012linear, nourian2013mm,bensoussan2016mean} in the LQ setting, \citet{sen2016mean,carmona2016probabilistic,carmona2017alternative} for the probabilistic analysis, \citet{lasry2018mean,cardaliaguet2020remarks} for the analytical viewpoint, \citet{carmona2021mean,dayanikli2023multi} for the incorporation of multiple major players, and~\citet{cui2024learning} for a learning algorithm). If the major player has a stochastic evolution, this may lead to a form of common noise for the minor players, which makes such models challenging to analyze and to solve. Stackelberg MFGs considers that a principal influences the agents in the mean field population (see e.g.~\citet{bensoussan2016mean,elie2019tale} for the theoretical background and~\citet{epidemics_SMFG,dayanikli2023machine} for machine learning methods).
The main difference between major-minor MFG and Stackelberg MFG is the equilibrium notion between the major agent/principal and the mean field population. In major-minor MFG, the major player and the minor players that constitute the large population are in a Nash equilibrium. However, in Stackelberg MFG, the principal and the players in the mean field population are in a Stackelberg equilibrium. Because of this difference, one of these notions can be a better fit for the application of interest. For example, if we are interested in public policy making, a Stackelberg MFG can be a better fit to model the interactions of a government and the individuals in a large society. Furthermore, other notions of solutions have been investigated, such as social optimum, also called mean field control (MFC) or control of McKean-Vlasov dynamics (see e.g.~\citet{bensoussan2013mean,carmona2013control} for the definitions, and~\citet{carmona2021convergence,carmona2022convergence,dayanikli2023deep} for deep learning methods), and correlated equilibrium (see~\citet{campi2022correlated,bonesini2024correlated} for theoretical background and~\citet{muller2022learningPSRO,muller2022learning} for learning algorithms). While Nash equilibrium assumes that the agents are in a non-cooperative setup, mean field control considers a purely cooperative setup and correlated equilibrium consider non-cooperative players who can share a common signal. These setups can be relevant depending on the application under consideration. For instance, MFC can be relevant in distributed robotics, since all the robots cooperate to achieve a common goal.

\section{Methods for discrete time models} 
\label{sec:methods-discrete-time}

In this section, we describe several methods to solve discrete time games. To fix the ideas, we present the algorithms in the context of mean field games, although similar ideas could be used for finite-player games. 
We start with a general description of fixed point methods to compute Nash equilibria. We then discuss two classical families of algorithms to compute best responses: policy optimization and dynamic programming.

\subsection{Fixed point algorithms}
\label{subsec:discrete_fixedpoint}

{\bf High-level description. }
Recall that the solution notion, namely Nash equilibrium, can be formulated as a fixed point problem. Therefore, the simplest idea is to iteratively apply the \textit{best response operator} and the \textit{mean field operator}. To be more specific, let us consider that, at the beginning of iteration $k \ge 0$, we have a candidate $\mu^k = (\mu^k_n)_{n=0,\dots,N_T}$ for the mean field sequence. We first compute a new policy, say $\pi^{k+1}$, which is a best response against $\mu^k$. Then we compute the new mean field, $\mu^{k+1}$ as the sequence induced when the whole population uses policy $\pi^{k+1}$. 
There are several methods to perform these steps in practice, depending on the setting. Specific methods will be presented and discussed in detail below, with a special attention given to methods allowing to compute a best response. 
This is summarized in Algorithm~\ref{algo:fixedp-iterations}. 

\begin{algorithm}[H]
\caption{Fixed point iterations \label{algo:fixedp-iterations}}

\textbf{Input:} initial mean field $\mu^0 = (\mu^0_n)_n$; number of iterations $K$

\textbf{Output:} equilibrium policy and mean field
\begin{algorithmic}[1]
\State{Compute a best response $\pi^0$ against mean field $\mu^0$}
\For{$k = 0,\dots, K-1$}
    \State{
    {\bf Mean field update: } Let $\mu^{k+1}$ be the mean field induced by policy $\pi^k$}
    \State{
    {\bf Policy update: } Compute a best response $\pi^{k+1}$ against mean field $\mu^{k+1}$}
\EndFor\\
\Return{$\pi^K, \mu^K$}
\end{algorithmic}
\end{algorithm}

\noindent{\bf Convergence. } A first remark is that, if the sequences of $(\mu^k)_{k \ge 0}$ and $(\pi^k)_{k \ge 0}$ converge, say respectively to $\mu$ and $\pi$, then $(\mu,\pi)$ form a Nash equilibrium. Indeed, in the limit, the mean field update does not modify $\mu$, which means that $\mu$ is the mean field generated by policy $\pi$; moreover, if the policy update does not modify $\pi$, it means that $\pi$ is a best response against $\mu$. We thus obtain a pair which satisfies the requirements of a Nash equilibrium, see Definitions~\ref{def:MFGNE-discrete} and~\ref{def:MFGNE-continuous}. This property is true even if there are multiple Nash equilibria.

A sufficient condition for the convergence of $\mu^k$ is that the composition $MF \circ BR$ of the two functions is a strict contraction. Intuitively, this holds provided each function is a strict contraction. In particular, the best response function needs to be Lipschitz in the mean field, with a small enough Lipschitz constant, which can be quite restrictive. One way to help this condition to be satisfied is to modify the reward by adding an extra entropic penalty, which regularizes the best response. Another option is to keep the original reward but to optimize only over a class of regularized policies such as policies that are soft-max of the Q-function. 
However, it should be noted that these two approaches modify the reward function and hence change the Nash equilibrium of the game, see e.g.~\citet{guo2019learning,cui2021approximately,guo2022entropy,anahtarci2023q}.   

As for the convergence of the policy, the main roadblock is that in many cases, there are multiple best responses to a given mean field. Here again, adding an entropic regularization in the reward or directly optimizing over smooth policies can help ensuring uniqueness of the best response at the expense of changing the Nash equilibrium.

Since computing a best response against a given mean field $\mu = (\mu_n)_{n=0,\dots,N_T}$ boils down to solving the MDP:
\begin{equation*}
    \max_{\pi} J(\pi;\mu).
\end{equation*}
\begin{remark}
    We assume that the maximum (or, if the objective is minimization, the minimum) of the objective is attained. Otherwise, we can aim to find the \textit{supremum} (or \textit{infimum}) and approximately optimal policies.
\end{remark}

This step can be tackled using any existing method for MDPs. Next, we describe two families of methods.

\noindent{\bf Variants. } 
Several variants of the above fixed point iterations have been introduced. In the discrete time setting, the two main methods are Fictitious Play and Online Mirror Descent. They were first studied in the continuous time and space setting e.g. by~\citet{cardaliaguet2017learning,hadikhanloo2017learning}.
\begin{itemize}
    \item {\bf Fictitious Play} consists in computing a BR not against the last mean field $\mu^k$ but against the average $\bar\mu^k$ of mean fields seen in past iterations, i.e., $\bar\mu^k_n(s) = \frac{1}{k} \sum_{j=1}^k \mu^j_n(s)$ for all $s,n$; see~\citet{elie2020convergence,perrin2020fictitious} for proofs of convergence under monotonicity condition and numerical experiments.
    \item {\bf Online Mirror Descent} consists in replacing the BR computation by a policy improvement based on first evaluating the current policy, computing the sum of all past value functions, and then letting the new policy be the softmax of this cumulative value function. To be specific, given the state-action value function $Q$ obtained at a given iteration (by evaluating the previous policy), we compute the next policy by taking $\pi(a|x) = (\mathrm{softmax}_\tau Q(x,\cdot))(a) = \exp(Q(x,a)/\tau) / \sum_{a' \in A} \exp(Q(x,a')/\tau)$, where $\tau>0$ is the temperature parameter. We refer to~\citet{perolat2022scaling} for more details and a proof of convergence under monotonicity assumption. 
\end{itemize}

\subsection{Computing a best response by policy optimization}

A first and somewhat straightforward approach to compute an optimal policy is to view the MDP as an optimization problem over the space of policies. If the state space or the action space is large, then it is common to replace the policy $\pi$ by a neural network $\pi_\theta$ with parameters $\theta \in \Theta$. Then, the problem becomes:
\begin{equation}
\label{eq:br-mdp-pitheta-mu}
    \max_{\theta \in \Theta} J(\pi_\theta;\mu).
\end{equation}
To find an optimal parameter, we can then use gradient-based methods for instance. In general, there is no closed-form formula for $\nabla_\theta J(\pi_\theta;\mu)$, and we need to estimate this term. Since $J$ is an expectation, we can use Monte Carlo samples and stochastic gradient descent (SGD) or extensions such as Adam optimizer. This is summarized in Algorithm~\ref{algo:sgd-policy-optim}. To alleviate the presentation, we introduce the following notation: a trajectory is denoted by $\tau = (s_n,a_n,r_n)_{n=0,\dots,N_T}$ with 
\begin{equation}
    \label{eq:discrete-time-MF-traj-simu}
    s_0 \sim \mu_0, \ a_n \sim \pi_{\theta,n}(s_n),\ s_{n}\sim P_{n-1}(s_{n-1}, a_{n-1}, \mu_{n-1}),\ r_n = r_n(s_n, a_n, \mu_n), n = 1,\dots,N_T.
\end{equation}
The total reward obtained for this trajectory is denoted by: 
\begin{equation}
    \label{eq:discrete-time-MF-total-reward-simu}
    \hat{J}(\tau) = \sum_{n=0}^{N_T} r_n.
\end{equation}
Note that $\tau$ implicitly depends on the policy that is used to generate the actions. Changing the policy implies changing the distribution of actions. To clarify the gradient's computation, one can introduce random variables $\epsilon_n$ and $U_n$ for the randomness coming from the dynamics and the policy and write:
\begin{equation}
\label{eq:discrete-time-state-action-sampling}
    s_{n} =  P_{n-1}(s_{n-1}, a_{n-1}, \mu_{n-1}, \epsilon_n),
    \qquad 
    a_n = \pi_{\theta,n}(s_n, U_n)
\end{equation}
where now $P$ and $\pi_{\theta,n}$ are deterministic functions. Now, we can compute the gradient of $\hat{J}$ with respect to  the policy's parameter $\theta$ for a given randomness realization $(\epsilon_n,U_n)_n$.

\begin{algorithm}[H]
\caption{BR computation by policy optimization \label{algo:sgd-policy-optim}}

\textbf{Input:} mean field $\mu = (\mu_n)_n$; number of iterations $K$; learning rates $(\beta^k)_{k=0,\dots,K-1}$

\textbf{Output:} (approximation of an) optimal neural network policy for~\eqref{eq:br-mdp-pitheta-mu}
\begin{algorithmic}[1]
\State Initialize $\theta^0$
\For{$k =0,\dots, K-1$}
    \State{Generate a trajectory: sample $(\epsilon_n)_n, (U_n)_n$ and use~\eqref{eq:discrete-time-state-action-sampling} to generate $\tau^k = (s_n,a_n,r_n)_{n=0,\dots,N_T}$ using~\eqref{eq:discrete-time-MF-traj-simu} with policy $\pi_{\theta^k}$ and mean field $\mu$}
    \State{Compute the gradient $\nabla_\theta \hat{J}(\tau^k)$, where $\hat{J}$ is defined in~\eqref{eq:discrete-time-MF-total-reward-simu}}
    \State{Update parameters: $\theta^{k+1} \leftarrow \theta^{k} + \beta^k \hat{J}(\tau)$}
\EndFor\\
\Return{$\theta^K$}
\end{algorithmic}
\end{algorithm}

The above policy optimization method (implicitly) relies on the knowledge of the model in order to compute the gradient $\nabla_\theta \hat{J}(\tau)$ because to compute this quantity exactly, one needs the reward and transition functions. However, the method can be adapted to model-free setting using policy gradient type methods. This class of methods builds upon the policy gradient theorem. We refer to the book of~\citet[Chapter 13]{sutton2018reinforcement}. Although we focused here on the MFG setting, RL policy optimization methods can also be applied in finite player games, see e.g.~\citet{hambly2023policy}. To ensure convergence of the algorithm, one needs to let the learning rates $(\beta^k)_k$ go to $0$. However, they cannot vanish too quickly, otherwise the learning stops too early. From the theoretical perspective, rigorous proofs of convergence can be obtained through stochastic approximation, which formulates precise conditions on the learning rates. We refer to the original work of~\citet{MR0042668}, and to the book of e.g.~\citet[Section 10.2]{MR2442439} for an application of stochastic approximation arguments to the convergence of stochastic gradient schemes.

\subsection{Computing a best response by dynamic programming}
\label{sec:BR-by-DP}
Another family of approaches consists in using the time structure of the MDP and using dynamic programming instead of viewing it as a pure optimization problem. To this end, we define the (state-only) value function as:
\[
    V^*_n(s;\mu) =  \max_{\pi} \ \EE\left[\sum_{\nprime=n}^{N_T} r_\nprime(s_\nprime, a_\nprime,\mu_\nprime) \,\big|\, s_n = s \right],
\]
which is the maximum possible expected total reward when the state trajectory starts at time $n$ from state $s$ given the mean field $\mu$. In particular, $\EE_{s \sim \mu_0}[V^*_0(s;\mu)]$ is the best possible value achievable by an individual player given the mean field $\mu$.

It can be shown that $V^* = (V^*_n)_n$ satisfies the Bellman equation, also called dynamic programming (DP) equation: 
\begin{empheq}[left = \empheqlbrace]{align}
    V^*_{N_T}(s;\mu) &= \max_{a \in A}\ r_{N_T}(s,a,\mu_{N_T})
    \\
    V^*_{n}(s;\mu) &= \max_{a \in A}\ \left\{ r_{n}(s,a,\mu_{n}) + \EE_{s_\nprime \sim P_n(s,a,\mu_n)}[V^*_{n+1}(s_\nprime;\mu)] \right\}.
\end{empheq}
Any policy whose support is over the set of maximizers in the above equation for each $s$ and $n$ is a best response against the mean field $\mu$. In fact, $\pi$ is a best response against $\mu$ if and only if: for every $n=0,\dots,N_T$, $s \in S$, $a \in A$,  
\[
    \pi_n(a|s,\mu_n) > 0 \Rightarrow a \in \argmax_{a \in A} \left\{ r_{n}(s,a,\mu_{n}) + \EE_{s_\nprime \sim P_n(s,a,\mu_n)}[V^*_{n+1}(s_\nprime;\mu)] \right\}.
\]
In particular, 
\[
    \EE_{s_0 \sim \mu_0}[V^*_0(s_0; \mu)] = \max_\pi J(\pi; \mu).
\]
In fact, it is often more convenient to work with the state-action value function, also called $Q$-function, which is defined as:
\[
    Q^*_n(s,a;\mu) =  \max_{\pi} \ \EE\left[\sum_{\nprime=n}^{N_T} r_\nprime(s_\nprime, a_\nprime,\mu_\nprime) \,\big|\, s_n = s, a_n = a \right].
\]
It corresponds to the maximum possible expected total reward when the state trajectory starts at time $n$ from state $s$ and the first action is $a$ given the mean field $\mu$, after which the player picks actions according to an optimal policy. We have the relation:
\[
    V^*_n(s;\mu) = \max_{\nu \in \mathcal{P}(A)} \EE_{a \sim \nu}[Q^*_n(s,a;\mu)],
\]
and $Q^*$ satisfies the Bellman equation (or DP equation):
\begin{empheq}[left = \empheqlbrace]{align}
    Q^*_{N_T}(s,a;\mu) &= r_{N_T}(s,a,\mu_{N_T})
    \label{eq:Qstar-terminal}
    \\
    Q^*_{n}(s,a;\mu) &= r_{n}(s,a,\mu_{n}) + \EE_{s_\nprime \sim P_n(s,a,\mu_n)}[\max_{a^\prime \in A}\  Q^*_{n+1}(s_\nprime,a^\prime;\mu)] .
    \label{eq:Qstar-n}
\end{empheq}
The main advantage of working with the $Q$-function is that, once $Q^*$ is computed, an optimal policy can easily be deduced: $\pi$ is an optimal policy provided it is distributed over the set of maximizers of $Q^*$. Namely, for every $n=0,\dots,N_T$, $s \in S$, $a \in A$, 
\[
    \pi_n(a|s) > 0 \Rightarrow a \in \argmax Q_n(s,\cdot;\mu).
\]
In fact, if one computes $Q^*$ by using the dynamic programming equation above, an optimal policy can be computed along the way. This is illustrated in Algorithm~\ref{algo:dynprog-Q}. Note that, in general, there are infinitely many best responses but perhaps only one of them gives the Nash equilibrium. Hence returning the $Q$-function instead of just one policy (chosen arbitrarily) can be useful.

\begin{algorithm}[H]
\caption{BR computation by dynamic programming\label{algo:dynprog-Q}}

\textbf{Input:} mean field $\mu = (\mu_n)_n$

\textbf{Output:} optimal $Q$-function and an optimal policy for~\eqref{eq:br-mdp-pitheta-mu}
\begin{algorithmic}[1]
\State{Let $Q^*_{N_T}(s,a;\mu)$ for all $s,a$ using~\eqref{eq:Qstar-terminal}}
\For{$n = N_T-1,\dots, 0$}
    \State{Compute $Q^*_{n}(s,a;\mu)$ for all $s,a$ using~\eqref{eq:Qstar-n}}
    \State{Let $\pi^*_n(\cdot|s)$ be an action distribution with support in $\argmax Q_n(s,\cdot;\mu)$ for all $s$}
\EndFor\\
\Return{$Q^*,\pi^*$}
\end{algorithmic}
\end{algorithm}

In Algorithm~\ref{algo:dynprog-Q}, each step is computed \textit{exactly} in the sense that we assume we can compute perfectly the expectation and the maximal value appearing in~\eqref{eq:Qstar-n}. In practice, if the state and action spaces are large, we can replace $Q^*$ and $\pi^*$ by neural networks and train them using using samples and an SGD-based method.

The above approach can also be used as a foundation for model-free algorithms. Computing the optimal $Q$-function can be done using the celebrated $Q$-learning algorithm~\citet{watkins1992q} (with tabular representation) or one of the more recent developments based on deep neural networks such as DQN~\citet{mnih2015human}. It is also possible to learn both a $Q$-function approximation and policy approximation using methods such as actor-critic algorithms. We refer to e.g.~\citet{arulkumaran2017deep,franccois2018introduction} for surveys on the basics of deep RL and we omit more recent developments for the sake of brevity.

\subsection{Numerical illustrations}
We now present two concrete examples that are solved numerically using learning methods. 
\subsubsection{Traffic routing with online mirror descent}
\label{sec:numerics-routing}
We present an illustration of the traffic routing model introduced by~\citet{cabannes2022solving} and discussed in Section~\ref{sec:example traffic}. We detail the model here. The network of roads is represented by an oriented graph. $L$ denotes the set of edges or links. The model presented by~\citet{cabannes2022solving} is in continuous time, but the implementation in OpenSpiel implementation uses a discrete time version of the model. Here, we directly focus on the discrete time version, in line with the notations introduced above. We fix a time horizon $T$ and an integer $N_T$, and we define the time step size $\Delta t = T/N_T$.  
\begin{itemize}
    \item Agents: cars (or drivers).
    \item States: the state of an individual at time step $n$ is $s_n = (\ell_n, w_n)$ corresponding to the current link and the waiting time. It takes  values in $X =  E\times \{0,\dots,N_T\}$.
    \item Actions: an individual can choose which edge to use next, among the set of successors of the current edge $\ell_n$. So the global action space is $A = E$, but there is a constraint at time step $n$ ensuring that the legal actions are only the successors of $\ell_n$. 
    \item Transitions: the initial state of a representative player is $(\ell_{0}, w_{0})$, and the dynamics are as follows: as long as $w_{n} > \Delta t$, the player cannot move to another link and the waiting time decreases, so:
    \[
        \ell_{n+1} = \ell_{n}, \qquad w_{n+1} = w_{n} - \Delta t,
    \]
    and when $w_{n} \le \Delta t$, the player can move to the next link and the waiting time is set to a new value, which depends on how congested the link is:
    \[
        \ell_{n+1} \sim \pi_{n}(s_{n}), \qquad w_{n+1} = c\left(\ell_{n},\mu_{n}(\ell_{n})\right),
    \]
    where $pi$ denotes the policy, $\mu_{n}(\ell_{n})$ is the total mass in link $\ell_n$ at time $n$, and $c: \mathcal{E} \times [0,1] \to \RR_+$ is a congestion function; the value $c\left(\ell, m\right)$ determines the waiting time when arriving in a link $\ell$ with a proportion $m$ of the total population. We refer to~\citet{cabannes2022solving} for more discussions on realistic choices of congestion functions $c$ based on empirical studies.
    \item Rewards: each player wants to arrive at the destination as soon as possible. Here, we assume that there is a common destination $D$ for all the drivers.\footnote{An extra component can be added to keep track of the destination if different groups of players have different destinations; see~\citet{cabannes2022solving} for more details.} So the reward is simply the total time spent on the road before reaching the destination. Mathematically, we take $r((\ell,w),a,\mu) = - \indic_{\ell \neq D}$, which makes the agent pay a cost of $1$ as long as the state is not the destination state (notice that here we view links and not nodes as states). We could also multiply by $\Delta t$ to take into account the time discretization but this will have no influence on the equilibrium policies. 
\end{itemize}

We solve the MFG by using the OMD algorithm mentioned in Section~\ref{subsec:discrete_fixedpoint}. We employ the implementation available in the open-source library OpenSpiel; see~\citet{lanctot2019openspiel}.\footnote{More precisely, we used the codes available at the following address (as of May 2024), which allows to reproduce the experiments of~\citet{cabannes2022solving}: \url{https://github.com/google-deepmind/open_spiel/tree/master/open_spiel/data/paper_data/routing_game_experiments}.} In particular, at each iteration, the value function and the mean field are computed using the transition matrices so there is no approximation. 

\begin{figure}
    \centering
    \includegraphics[width=0.3\linewidth]{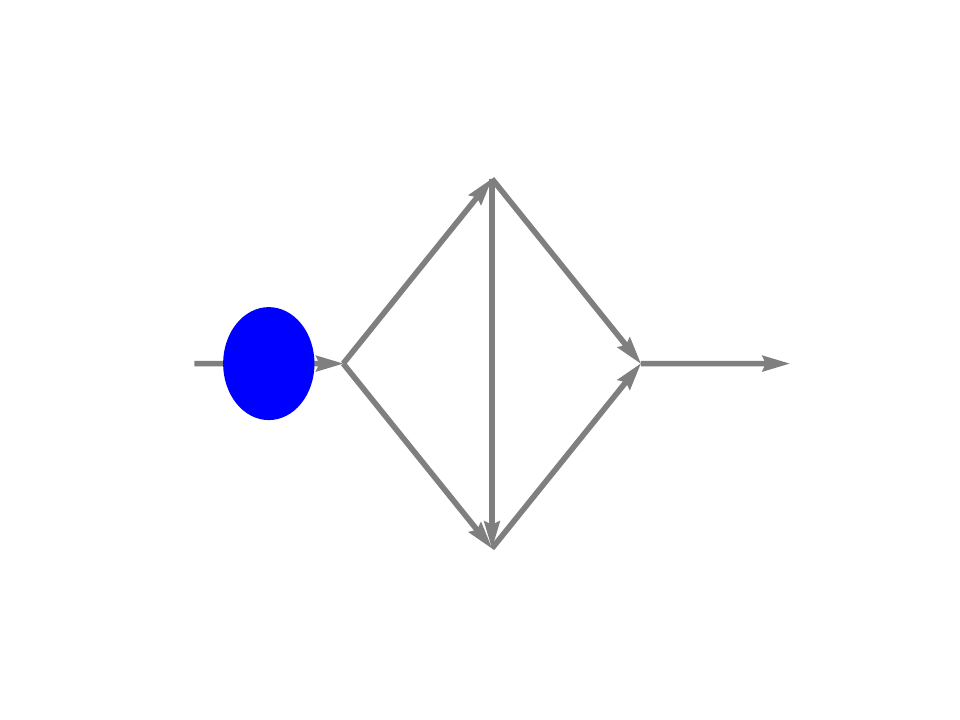}
    \includegraphics[width=0.3\linewidth]{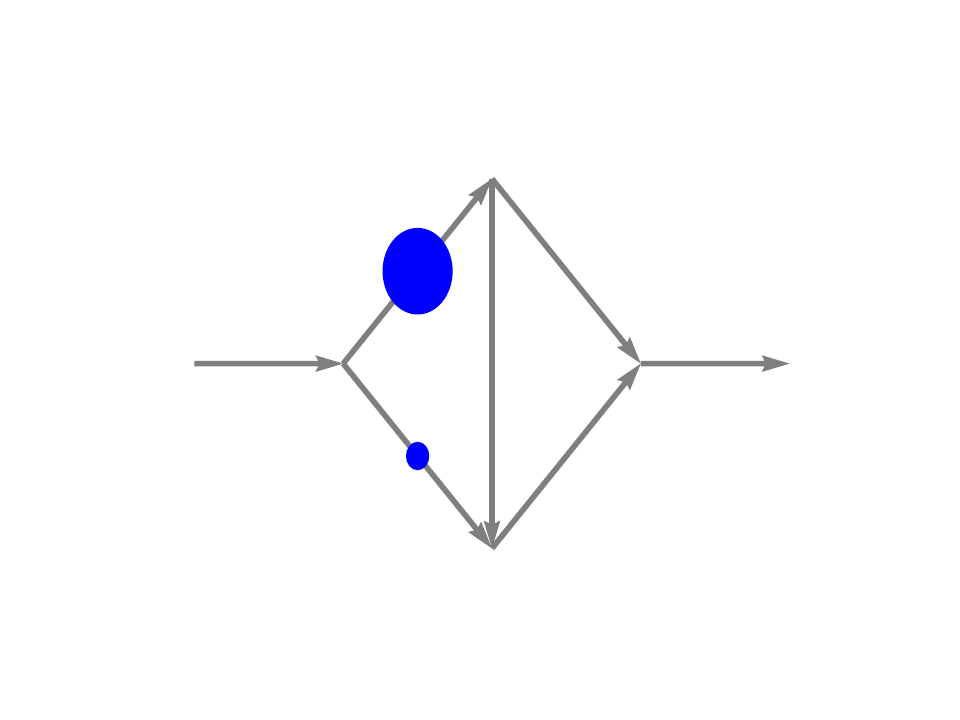}
    \includegraphics[width=0.3\linewidth]{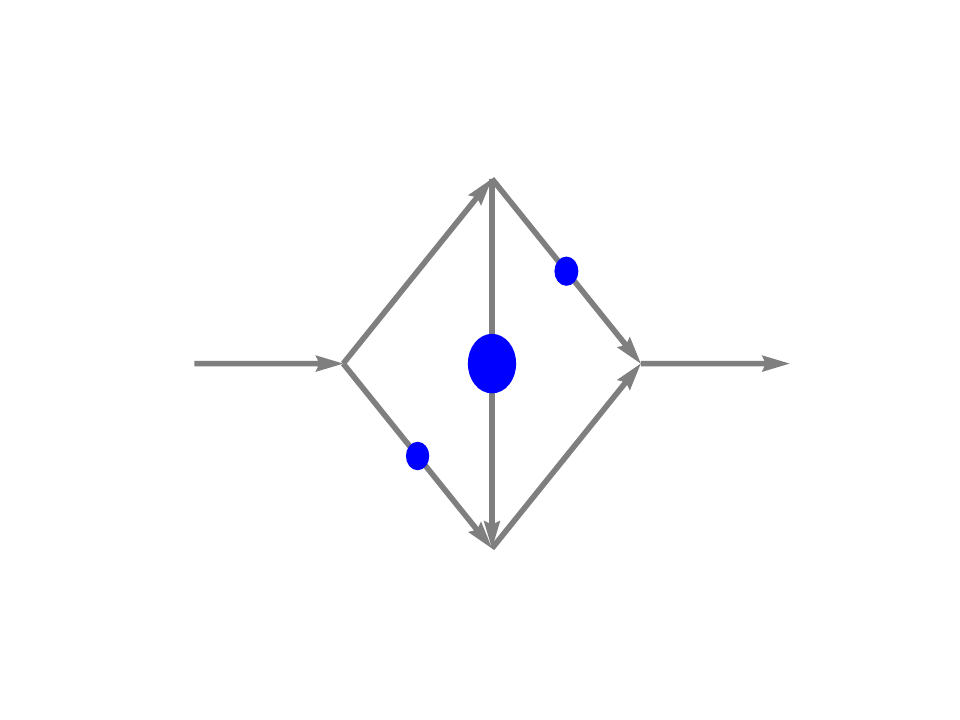}
    \includegraphics[width=0.3\linewidth]{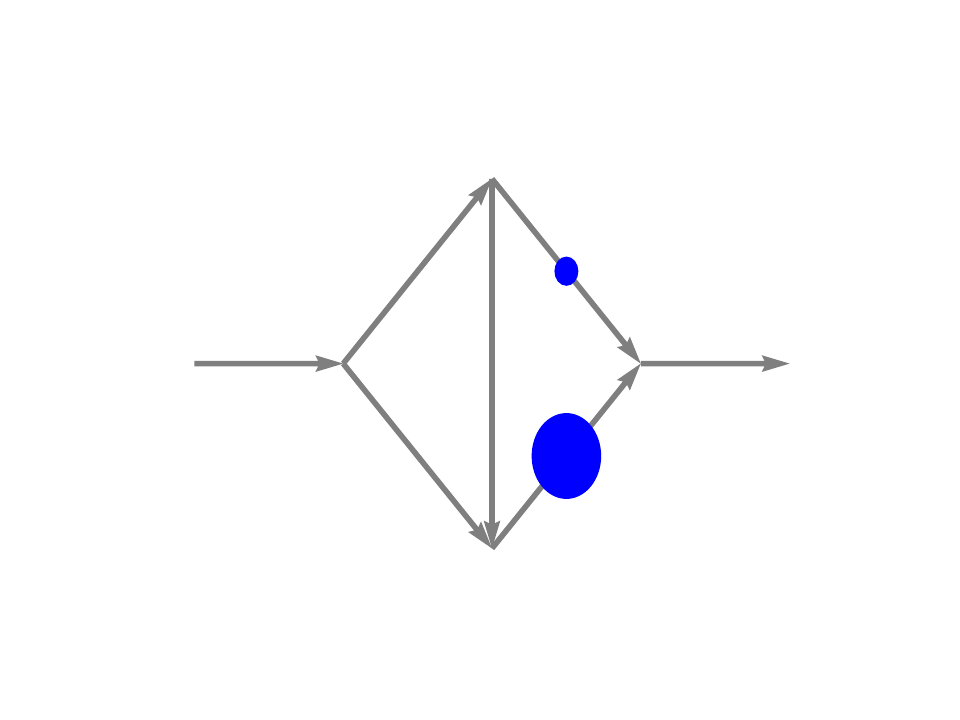}
    \includegraphics[width=0.3\linewidth]{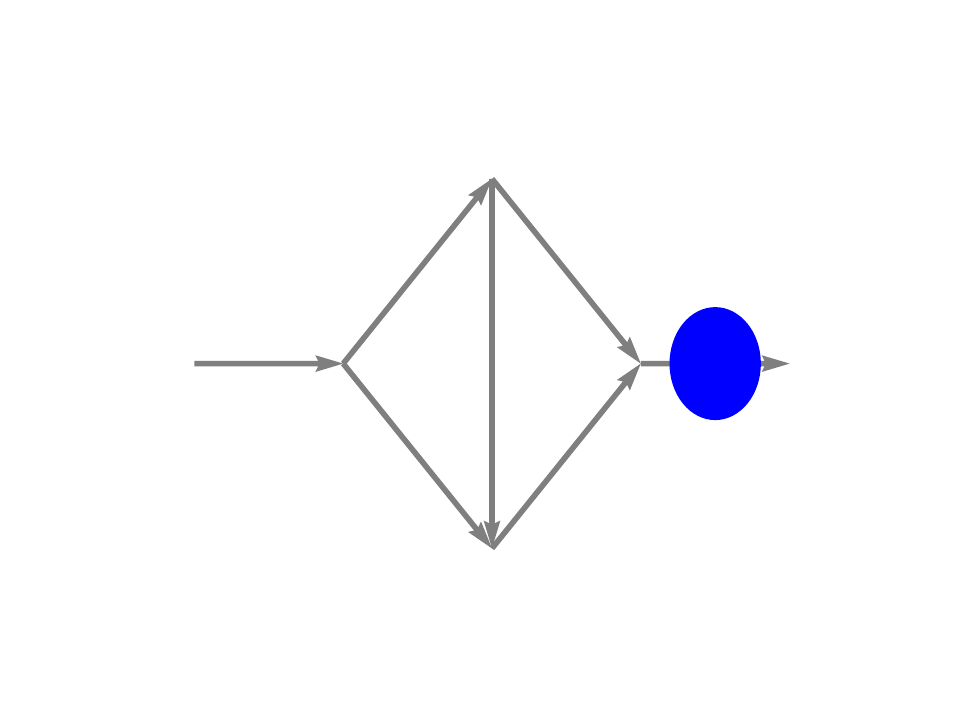}
    \caption{Evolution of the population distribution in the MFG of Braess network, at times $n=0,1,8,10,16$ (from top left to bottom right). The size of the disc represents the proportion of the population on the correponsding link at that time. }
    \label{fig:braess_dynamics}
\end{figure}

As for the network, we consider here an MFG version of the celebrated Braess paradox~\citet{braess2005paradox}. The network and the evolution of the population using the equilibrium policy learnt by OMD are shown in Figure~\ref{fig:braess_dynamics}. In this example, there are 20 time steps, a single origin and a single destination. We refer to~\citet{cabannes2022solving} for more details. 

\begin{figure}
    \centering
    \includegraphics[width=0.6\linewidth]{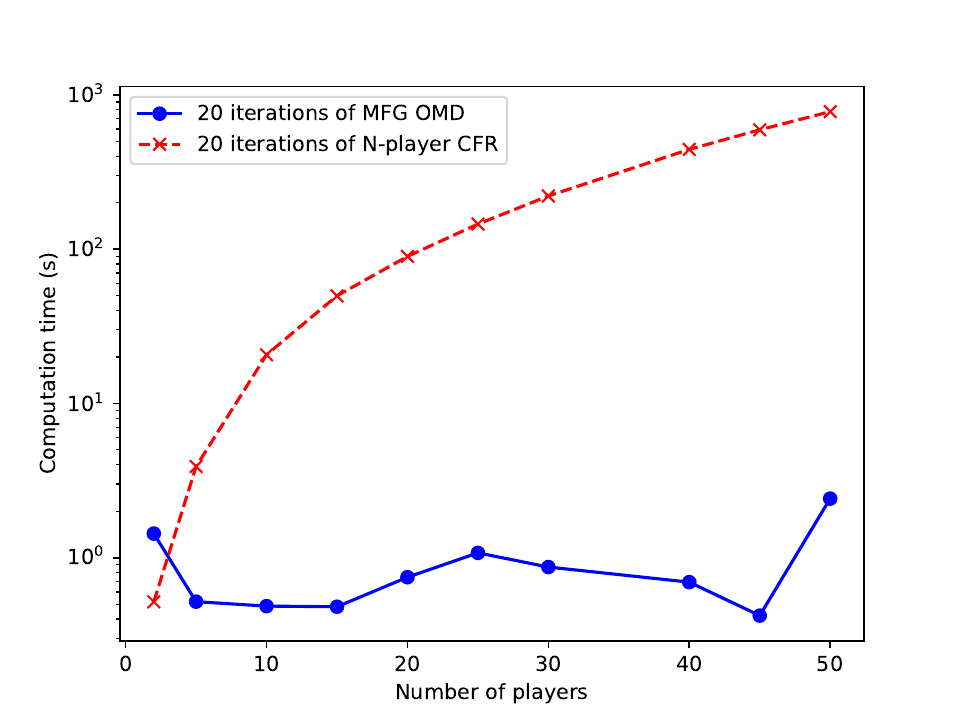}
    \caption{Computational time of CFR for the $N$-player Braess game versus OMD for the MFG. The results are average over $10$ runs.}
    \label{fig:computation_time_algo_n_vs_mean_field}
\end{figure}

Next, we illustrate the advantage of solving an MFG instead of a finite player game, even when one is interested in a scenario with a finite population. For various values of $N$, we compare the computational time of two approaches: (1) solving the $N$-player Braess network game using counterfactual regret minimization with external sampling (ext CFR), and (2) solving the MFG using OMD. For each $N$, we run the computation $10$ times and report the average. The results are reported in Figure~\ref{fig:computation_time_algo_n_vs_mean_field}. The computational time of CFR increases exponentially with $N$, while the cost of solving the MFG is independent of $N$.

\subsubsection{Crowd motion with fictitious play and RL} 
\label{sec:numerics-crowd}
In this section, we solve an MFG modeling crowd motion in a grid world by using fictitious play (FP) described in Section~\ref{subsec:discrete_fixedpoint} and reinforcement learning. We consider the model discussed in Section~\ref{sec:ex-crowd} with a 2D grid world and 5 possible actions. The noise takes value $(0,0)$ with probability $1-p$ and takes value any of the 4 directions with probability $p/4$. The reward is simply: $-\mu(x)$, which reflects crowd aversion. The domain is represents 4 rooms connected by doors. This model was introduced by~\citet{geist2022concave}, who studied the connection with RL. It has then be used as a benchmark problem in several other works, e.g.~\citet{algumaei2023regularization}. In the simulations, we consider a grid world of size $9 \times 9$; we took $p=0.2$ for the noise probability and horizon $N_T=15$.

We solve the problem using three different algorithms based on fictitious play. In each case, the mean field evolution is computed exactly using the transition matrix, but the methods differ by the way the best response (BR) and the average policy are computed:
\begin{itemize}
    \item Exact Fictitious Play: standard fictitious play, with perfect computation of the BR (using backward induction);
    \item Tabular Q-learning Fictitious Play: computation of the BR using tabular Q learning and computation of the average policy by a tabular representation;
    \item Deep Fictitious Play: computation of the BR using deep RL (DQN) and computation of the average policy by deep learning (see~\citet{lauriere2022scalable} for more details).
\end{itemize}  

\begin{figure}
    \centering
    \includegraphics[width=0.6\linewidth]{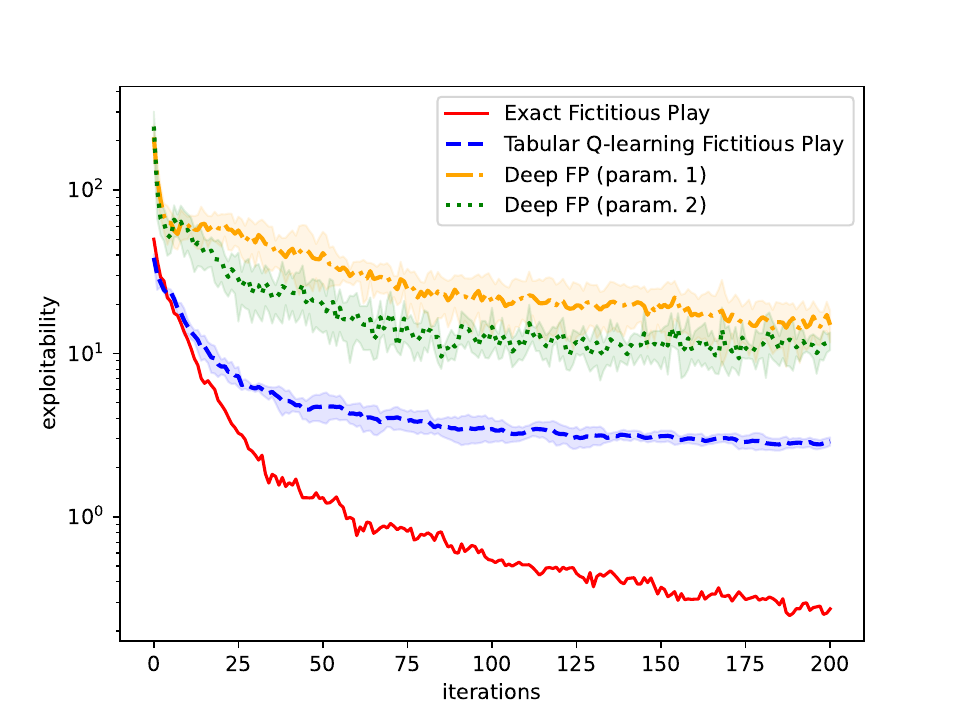}
    \caption{Exploitabilities for the 4-room crowd motion example with several variants of fictitious play, as described in the text.}
    \label{fig:crowd-openspiel-fp-rl-exploitabilities}
\end{figure}

The simulations are done using the library OpenSpiel (see~\citet{lanctot2019openspiel}), which contains the MFG crowd motion model as well as the above algorithms.\footnote{See \url{https://github.com/google-deepmind/open_spiel/tree/master/open_spiel/python/mfg}.}
Figure~\ref{fig:crowd-openspiel-fp-rl-exploitabilities} shows the exploitabilities. For the methods based on RL, we ran the algorithms 10 times and computed the average (lines) as well as the standard deviation (shaded areas). For the tabular Q learning method, we used the default parameters in OpenSpiel. For the deep RL method, we ran it with two different choices of hyperparameters: in the first, we use $50$ inner iterations of DQN, an architecture of $[64,64]$ neurons and a constant $\epsilon=0.1$; in the second, we use $500$ inner iterations, $[128,128,128]$ neurons and an $\epsilon$ decaying from $0.1$ to $0.001$.\footnote{The code is available upon request.}
We observe that standard fictitious play converges faster than the methods where the BR is approximated using RL. Among the two methods using RL, the one with tabular Q-learning converges faster than the one using deep RL. This is probably due to the fact the game is small enough that tabular methods can perform very well. For larger games, we expect deep RL methods to outperform tabular ones. Among the two sets of hyperparameters, we observe -- as expected -- that the second set of hyperparameters leads to a better performance. In any case, these results are provided only for the sake of illustration; to obtain better results, one should tune the hyperparameters of the tabular RL and the deep RL methods.

\begin{figure}
    \centering
    \includegraphics[width=0.9\linewidth]{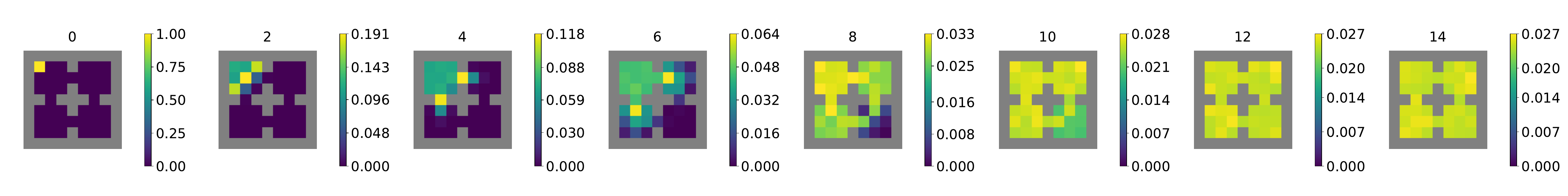}

    \includegraphics[width=0.9\linewidth]{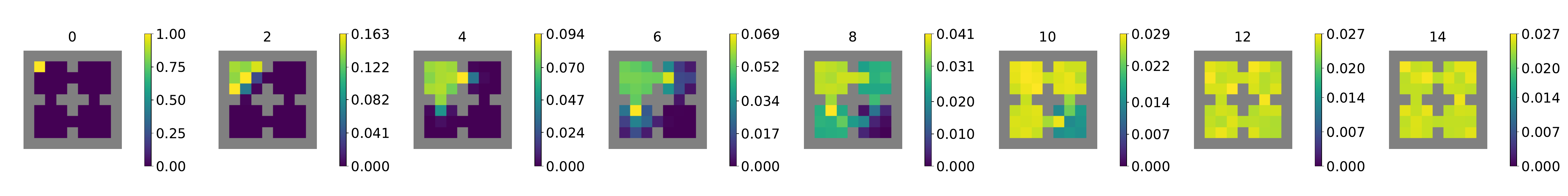}

    \includegraphics[width=0.9\linewidth]{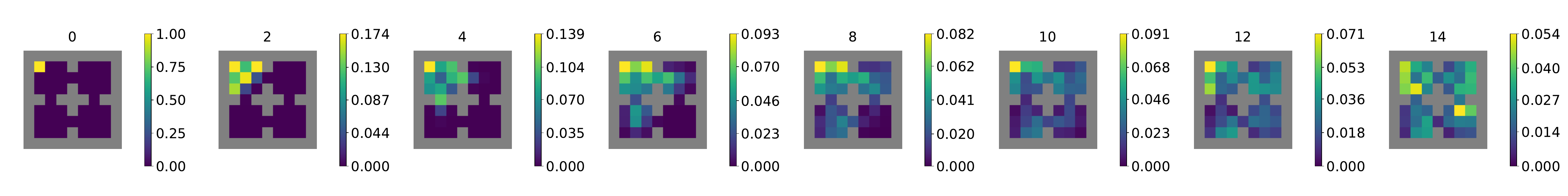}
    
    \includegraphics[width=0.9\linewidth]{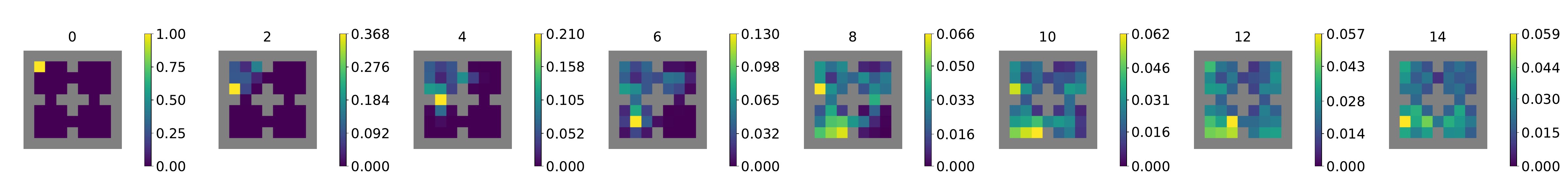}
    \caption{Evolution of the distribution for the 4-room crowd motion example, solved with exact FP (row 1), FP with tabular Q learning (row 2), deep FP with hyperparameters 1 (row 3) and deep FP with hyperparameters 2 (row 4). Each column corresponds to a different timestep.}
    \label{fig:crowd-openspiel-fp-rl-distributions}
\end{figure}

Figure~\ref{fig:crowd-openspiel-fp-rl-distributions} displays the evolution of the distribution when using the policy learnt by each of the algorithm (for the RL methods, we show one the result for the policy learnt in one of the 10 runs). We see that the agents spread through the four rooms in order to decrease the population density. The final distribution is more uniform with the exact and tabular RL methods and with the deep RL method (in the absence of hyperparameters tuning).

\section{Methods for continuous time models} 
\label{sec:methods-continuous-time}

We now turn our attention to machine learning methods to solve continuous time models as presented in Section~\ref{sec:finite-mf-games}. A first possibility is to discretize the time, and then use one of the methods introduced for discrete time games (see Section~\ref{sec:methods-discrete-time}). For instance, iterative methods combined with deep learning for finite player games have been proposed and analyzed by~\citet{hu2021deep,han2022convergence}. For the sake of brevity, we omit these methods and we focus on approaches that are really specific to continuous time models. 

An important advantage of using a continuous time formulation is the possibility to phrase optimality conditions through differential equations, and then solve numerically these equations. Below, we first discuss methods based on stochastic differential equations (SDEs) and then methods based on partial differential equations (PDEs). For simplicity of notation, we will take the volatility as a constant function, i.e., $\sigma(t, X_t, \mu_t)=\sigma$. The methods can be extended to more general forms of volatility.

\subsection{Deep learning methods for solving FBSDEs}
\label{subsec:deeplearning_fbsde}

Under certain technical assumptions, see e.g.~\citet[Chapter 3 and 4]{carmona2018probabilistic}, a solution of an MFG (i.e., mean field Nash equilibrium) can be characterized with forward-backward differential equations (FBSDEs) that are of the form:
\begin{equation}
\label{eq:fbsde_general}
\begin{aligned}
    dX_t &= b(t, X_t, \hat{\alpha}_t, \mu_t) dt + \sigma dW_t, && X_0=\zeta\sim\mu_0 \\
    dY_t &= - \partial_x H(t, X_t, \hat{\alpha}_t, \mu_t, Y_t)dt + Z_t dW_t,\quad && Y_T= \partial_x g(X_T, \mu_T),
\end{aligned}
\end{equation}
where $\mu_t= \mathcal{L}(X_t)$ and $\hat{\alpha}_t$ is the unique minimizer of the Hamiltonian, i.e., $$\hat{\alpha} = \argmin_{\alpha_t} H(t, X_t, \alpha_t, \mu_t, Y_t),$$
where the Hamiltonian is defined as
\begin{equation*}
    H(t, X_t, \alpha_t, \mu_t, Y_t) = b(t, X_t, \alpha_t, \mu_t) Y_t + f(t, X_t, \alpha_t, \mu_t).
\end{equation*}
This shows that $\hat{\alpha}_t$ can be expressed as a function of the state of the representative agent $X_t$, the population distribution $\mu_t$, and the \textit{adjoint} process $Y_t$: $\hat \alpha_t=\hat \alpha (t, X_t, \mu_t, Y_t)$. In this formulation the adjoint process $\bY=(Y_t)_{t\in[0,T]}$ characterizes the derivative of the value function $u: (t,x) \mapsto (t,x)$ with respect to the state variable, where the value function is defined as:
\begin{equation*}
\label{eq:value}
    u(t, x) = \min_{(\alpha_s)_{s\in[t, T]}} \EE\Big[\int_t^T f(s, X_s, \alpha_s, \mu_s)ds + g(X_T, \mu_T)\big| X_t=x \Big].
\end{equation*}
As mentioned in the introduction of this section, for the sake of simplicity in the presentation of the problem, we assume that the volatility is constant, i.e., in the form of $\sigma(t, X_t, \mu_t)=\sigma$. The derivation of the above FBSDE~\eqref{eq:fbsde_general} that characterizes the mean field Nash equilibrium follows from Pontryagin stochastic maximum principle. Their derivation and the analysis of the existence and uniqueness of their solutions are beyond the scope of this tutorial, which aims to focus on the machine learning methods to solve large population games. Interested readers can refer to~\citet[Chapter 3]{carmona2018probabilistic}. Therefore, we will directly focus on solving the FBSDE system with machine learning techniques after introducing the FBSDE model for the simple example given in Section~\ref{cont_ex:employee}.

\begin{remark}[Example] We recall that in the example given in Section~\ref{cont_ex:employee}, the running cost is $f(t, X_t, \alpha_t, \mu_t) = \frac{1}{2}\alpha_t^2-U(X_t)$ and the drift is $b(t, X_t, \alpha_t, \mu_t) = \alpha_t +\bar{X_t}$ where $\bar{X}_t=\int_{\RR} x\mu_t(x)$. The mean field Nash equilibrium control is then given as the minimizer of the Hamiltonian:
$
    \hat{\alpha}_t = -Y_t,
$ 
where $(\bX, \bY, \bZ)$ solves the following FBSDE: 
    \begin{equation*}
        \begin{aligned}
            dX_t &= (-Y_t + \bar{X}_t) dt + \sigma dW_t, \qquad&& X_0 =\zeta\sim \mu_0\\
            dY_t &= \partial_x U(X_t)dt + Z_t dW_t,\quad && Y_T=-\partial_x U(X_T).
        \end{aligned}
    \end{equation*}
\end{remark}

The FBSDE~\eqref{eq:fbsde_general} characterizing the mean field Nash equilibrium has two important challenges. First, there are forward and backward components which are coupled. The forward component moves forward in time starting from the initial time $t=0$ and the backward component moves backward in time starting from the terminal time $t=T$. Second, these are McKean-Vlasov type equations which means that the distribution $\bmu=\mathcal{L}(X_t)_{t\in[0,T]}$ of the process $\bX$ is in the dynamics. If the model is linear-quadratic, the corresponding FBSDE can be solved by proposing an \textit{ansatz} for $Y_t$, for example in the form of $Y_t=A_t X_t+B_t \bar{X}_t +C_t$ where $A_t, B_t,$ and $C_t$ are deterministic functions of time; see e.g.~\citet{malhame2020mean}. However, obtaining explicit solutions of FBSDEs is in general a difficult task and instead numerical approaches can be used. In the following subsections we will show how the FBSDEs that characterize the solution of MFGs can be solved by using machine learning.

\begin{remark}
    In this section, we focus on the FBSDE system~\eqref{eq:fbsde_general} that comes from Pontryagin stochastic maximum principle. However, we could have used also another probabilistic approach where the adjoint process characterizes the value function; see~\cite[Section 4.4]{carmona2018probabilistic}. This yields a different McKean-Vlasov FBSDE system. However, the numerical approaches introduced in this section can be adapted in a straightforward way to solve this FBSDE system.
\end{remark}

\subsubsection{Algorithm 1: Iterative learning}
\label{sec:FBSDE-iterative-algo}
The first approach we introduce to solve FBSDE~\eqref{eq:fbsde_general} with machine learning methods is an \textit{iterative} approach, which is based on a similar idea as the fixed point algorithms introduced in Section~\ref{subsec:discrete_fixedpoint}. In order to be able to implement the numerical approach, the first step is to discretize the time and to write the forward SDE in discrete time form. The second step is to fix a flow for backward process $\bY$ and to simulate a population of size $N$ in order to calculate the empirical distribution of the state $X_t$ as an approximation for $\mu_t$. The discretized forward equation can be written as follows:
\begin{equation}
\label{eq:discrete_FBSDE_forward}
    X^i_{t_{n+1}} = X^i_{t_n} + \Delta t \times b(t_n, X^i_{t_n}, \hat \alpha^i_{t_n}, \mu^N_{t_n})+\sigma \Delta W^i_{t_n}, \quad X^i_0 \sim \mu_0 
\end{equation}
for all $i\in [N]$ where $\hat \alpha^i_{t_n} = \hat{\alpha}(t_n, X^i_{t_n}, \mu^N_{t_n}, Y^i_{t_n})$ and $\Delta W^i_{t_n}\sim \mathcal{N}(0,\Delta t)$ with $\Delta t = t_{n+1}-t_n$ where $t_n = n \Delta t$ represents the discretized time steps in $\mathcal{T} = \{0, t_1, t_2, \dots, T\}$. The mean field is approximated by the empirical distribution of states of the Monte Carlo simulated population, i.e., $\mu_{t_n}^N = \frac{1}{N}\sum_{i=1}^N\delta_{X_{t_n}^i}$. We emphasize that since $\hat \alpha^i_{t_n}$ is a function of time $t_n,$ state $X^i_{t_n},$ mean field $\mu^N_{t_n},$ and adjoint process $Y^i_{t_n}$, the function $b$ can also be written as a function of $t_n, X^i_{t_n}, \mu^N_{t_n}, Y^i_{t_n}$. After the mean field is approximated, we need to \textit{learn} the solution of the backward equation. In order to do this, we follow the strategy proposed by the deep BSDE method~\citet{han2018solving} and we use a \textit{shooting} method. In order to introduce this idea, we write the BSDE in~\eqref{eq:fbsde_general} in integral form:
\begin{equation*}
    Y_t = \partial_x g(X_T) + \int_t^T \partial_x H(s, X_s, \hat{\alpha}_s, \mu_s)ds - \int_t^T Z_s dW_s,
\end{equation*}
where $\hat \alpha_{s} = \hat{\alpha}(s, X_{s}, \mu_{s}, Y_{s})$. In this backward equation, process starts from the terminal condition $Y_T = \partial_x g(X_T, \mu_T)$ at the terminal time $T$ and moves backward in time. Instead, we can start from an initial point $Y_0$ (which needs to be chosen, i.e., a control variable that we will optimize) and follow the same dynamics to arrive at a terminal point which we aim to match with the terminal condition. Therefore, we write the same differential equation forward in time by using the following integral form:
\begin{equation*}
    Y_t = Y_0 - \int_0^t \partial_x H(s, X_s, \hat{\alpha}_s, \mu_s)ds + \int_0^t Z_s dW_s.
\end{equation*}

Now, the objective is to \textit{shoot} for the terminal condition by minimizing the distance between the simulated terminal $Y_T$ and the terminal condition $\partial_x g(X_T, \mu_T)$ by choosing the initial point $Y_0$ and the coefficient $\bZ$ in front of the randomness. Then, the problem becomes (after discretizing the time and approximating the mean field with the empirical distribution $\mu^N_{t_n}$):
\begin{equation}
\label{eq:shooting_Y}
\begin{aligned}
    \min_{Y_0,  \bZ}\ &\EE\big[\big(Y_T-\partial_x g(X_T, \mu_T^{N})\big)^2\big]\\
    \text{s.t: } &Y_{t_{n+1}} = Y_{t_n} - \Delta t \times \partial_x H(t_n, X_{t_n}, \hat \alpha_{t_n}, \mu^N_{t_n})+Z_{t_n} \Delta W_{t_n},\ \forall n \in \mathcal{T} \backslash \{T\},
\end{aligned}
\end{equation}
given the Monte Carlo simulated trajectories $(X_{t_n}^i)_{t_n\in\mathcal{T}, i\in[N]}$ and $(\mu_{t_n}^N)_{t_n\in\mathcal{T}}$. For learning the solution of the problem~\eqref{eq:shooting_Y}, we will replace controls $(Y_0, \bZ)$ by parameterized functions $y_{0, \theta_1} : \RR \rightarrow \RR$ and $z_{\theta_2}: [0,T] \times \RR \rightarrow \RR$. In order to approximate these functions we will use neural networks. Here, the initial point $Y_0$ is implemented as a function of initial state $X_0$, since it represents the value at time $t=0$, which is a function of the state at time $0$. The process $\bZ$ is implemented as a function of time $t$ and the state $X_t$, since it is a stochastic process. In order to update the neural network parameters to minimize the cost in~\eqref{eq:shooting_Y}, we can use SGD or mini batch gradient descent. More precisely, the objective is to minimize over the neural network parameters $\theta = (\theta_1, \theta_2)$ the (empirical) objective:
\begin{equation}
\label{eq:averaged_shooting_Y}
\begin{aligned}
    \mathbb{J}^N(\theta) = \EE \left[\dfrac{1}{N}\sum_{i=1}^N \big(Y^{i,\theta}_T-\partial_x g(X^i_T, \mu_T^{N})\big)^2\right],
\end{aligned}
\end{equation}
given the previously Monte Carlo simulated trajectories $(X_{t_n}^i)_{t_n\in\mathcal{T}, i\in[N]}$ and $(\mu_{t_n}^N)_{t_n\in\mathcal{T}}$.

After learning $\bY=(Y^i_{t_n})_{t_n\in\mathcal{T}, i\in[N]}$ with the shooting method, we can update $\bX=(X_{t_n}^i)_{t_n\in\mathcal{T}, i\in[N]}$ by using~\eqref{eq:discrete_FBSDE_forward} given $\bY$. We can continue the iterations until convergence. The details of the pseudo code can be found in Algorithms~\ref{algo:simu-forward-XY_Y},~\ref{algo:SGD-MFG-Alg1}, and~\ref{algo:convergence-iterations}.

    The above algorithm iteratively compute an approximate best response and the induced approximate mean field. It can be expected to converge to yield an approximate Nash equilibrium under suitable contraction assumptions. It can be modified based on ideas similar to fictitious play or online mirror descent (see Section~\ref{subsec:discrete_fixedpoint}), although these are not always straigthforward to implement with deep networks. Another approach, which directly aims for the FBSDE solution is discussed below.

    \begin{algorithm}[H]
    \caption{Monte Carlo simulation of an interacting batch $\bY$ \label{algo:simu-forward-XY_Y}}
    
    \textbf{Input:} Number of particles $N$; state trajectories $\bX$ for $N$ particles, mean field trajectory $\bmu^{N}$; time horizon $T$; time increment $\Delta t$; initial distribution $\mu_{0}$; control functions $y_0,z$
    
    \textbf{Output:} Approximate sampled trajectories of $(\bY,\bZ)$.
    \begin{algorithmic}[1]
    \State{Let $n=0$, $t_0 = 0$ and set $Y^i_0 = y_0(X^i_0), \forall i \in [N]$}
    \While{$n \times \Delta t = t_n < T$}
        \State{Set $ Z^i_{t_n} = z({t_n}, X^i_{t_n}),\ i \in [N]$}\vskip1mm
        \State{Set $\hat \alpha^i_{t_n} = \hat{\alpha}({t_n},X^i_{t_n},{\mu}^N_{t_n}, Y^i_{t_n}),\ i \in [N]$}\vskip1mm
        \State{Sample $\Delta W^i_{t_n} \sim \mathcal{N}(0,\Delta t),\ i \in [N]$ }\vskip1mm
        \State{Let 
            $
                Y^{i}_{t_{n+1}} = Y^{i}_{t_n} - \partial_x H(t_n, X^{i}_{t_n}, \hat \alpha^{i}_{t_n}, \mu^N_{t_n})\Delta t + Z^{i}_{t_n} \Delta W^{i}_{t_n}
            $, $ i \in [N]$
        }
        \State{Set $t_n = t_n+\Delta t$ and $n=n+1$}
    \EndWhile
    \State{Set $N_T = n,\ t_{N_T} = T$}\\
    \Return $(Y^i_{t_n},Z^i_{t_n})_{n =0,\dots,N_T,\ i \in [N]}$ 
    \end{algorithmic}
    \end{algorithm}

\begin{algorithm}[h]
\caption{Stochastic Gradient Descent for solving~\eqref{eq:shooting_Y}\label{algo:SGD-MFG-Alg1}}

\textbf{Input:} Initial parameter $\theta_0$; number of iterations ${K}$; sequence $(\beta_{{k}})_{{k} = 0, \dots, {K}-1}$ of learning rates; number of particles $N$; state trajectories $\bX$ for $N$ particles, mean field trajectory $\bmu^{N}$; time horizon $T$; time increment $\Delta t$; initial distribution $\mu_0$

\textbf{Output:} Approximation of $\theta^*$ minimizing $\mathbb{J}^N$ defined by~\eqref{eq:averaged_shooting_Y}
\begin{algorithmic}[1]

\For{${k} = 0, 1, 2, \dots, {K}-1$}
    \State \parbox[t]{\dimexpr\textwidth-\leftmargin-\labelsep-\labelwidth}{Sample $S = (Y^i_{t_{{n}}},Z^i_{t_{{n}}})_{{n} =0,\dots,{N_T},\ i \in  [N]}$ 
    using Algorithm~\ref{algo:simu-forward-XY_Y} with control functions $(y_0, \boldsymbol z) = (y_{0,\theta_{{k},1}}, z_{\theta_{{k},2}})$ and parameters: number of particles $N$; state trajectories $\bX$ for $N$ particles, mean field trajectory $\bmu^{N}$; time horizon $T$; time increment $\Delta t$; initial distribution $\mu_{0}$ \strut}
    \State{Compute the gradient $\nabla \mathbb{J}^N(\theta_{{k}})$ of $\mathbb{J}^N(\theta_{{k}})$}
    \State{Set $\theta_{{k}+1} =  \theta_{{k}} -\beta_{{k}} \nabla \mathbb{J}^N(\theta_{{k}})$ }
\EndFor

\State \Return $\theta_{{K}}$ and trajectories $\bY^{\theta_K}, \bZ^{\theta_K}$ for $N$ particles.
\end{algorithmic}
\end{algorithm}

\begin{algorithm}[h]
\caption{Iterative Learning for MFG\label{algo:convergence-iterations}}

\textbf{Input:} Number of particles $N$; time horizon $T={t_{N_T}}$; time increment $\Delta t$; initial distribution $\mu_0$; number of iterations $K$

\textbf{Output:} Approximation of equilibrium mean field flow $\bmu^{*}$
\begin{algorithmic}[1]

\State{Initialize $\bY^{0}=(Y^{i,0}_{t_0}, \dots, Y^{i,0}_{t_{N_t}})_{i \in [N]}$ arbitrarily 
}
\For{$k=1,\dots,K$}
    \State \parbox[t]{\dimexpr\textwidth-\leftmargin-\labelsep-\labelwidth}{
    Pick $X_0^{i,k}\sim\mu_0$ i.i.d, set $n=0$\strut}\vskip1mm
    \While{$n \times \Delta t = t_n < T$}\vskip1mm
        \State{Let $ \mu^{N,k}_{t_n} = \frac{1}{N} \sum_{i=1}^N \delta_{X^{i,k}_{t_n}}$ }\vskip1mm
        \State{Set $\hat \alpha^{i,k}_{t_n} = \hat{\alpha}({t_n},X^{i,k}_{t_n},{\mu}^{N,k}_{t_n}, Y^{i,(k-1)}_{t_n}),\ i \in [N]$}\vskip1mm
        \State{Sample $\Delta W^i_{t_n} \sim \mathcal{N}(0,\Delta t),\ i \in [N]$ }\vskip1mm
        \State{Let 
            $
                X^{i,k}_{t_{n+1}} = X^{i,k}_{t_n} + b(t_n, X^{i,k}_{t_n}, \hat\alpha^{i,k}_{t_n}, \mu^{N,k}_{t_n})\Delta t + \sigma \Delta W^{i}_{t_n}
            $, $ i \in [N]$}\vskip1mm
        \State{Set $t_n = t_n+\Delta t$ and $n=n+1$}\vskip1mm
    \EndWhile
    \State{Learn trajectory of $\bY^{k}$ and get $\theta^k$ by using Algorithm~\ref{algo:SGD-MFG-Alg1} given $\bX^{k}$, and $\bmu^{N,k}$ }\vskip1mm 
\EndFor
\State Set $\bmu^{N,*} = \bmu^{N,k}$\vskip1mm
\State \Return $\bmu^{N,*}$ and $\theta^K$
\end{algorithmic}
\end{algorithm}

\subsubsection{Algorithm 2: Simultaneous learning}
\label{sec:fbsde-algo2-simultaneous}
The second approach we introduce to solve FBSDE~\eqref{eq:fbsde_general} with machine learning methods tries to overcome the possible inefficiency of the iterations by simulating the coupled processes $\bX$ and $\bY$ simultaneously. In FBSDE~\eqref{eq:fbsde_general} (and also in its discretized version), the backward and forward components are coupled; therefore, we cannot simulate them directly. In order to be able to do this, we need to rewrite the backward differential equation forward in time as introduced previously. The aim is to \textit{shoot} for the terminal condition of the backward equation as introduced in the previous algorithm. However, instead of alternating between the forward and the backward dynamics, now we simulate the dynamics of $\bX$ and $\bY$ at the same time. 

To be specific, the problem becomes:
\begin{equation}
\label{eq:shooting}
    \begin{aligned}
             \min_{Y_0, \bZ}\ & \EE\big[\big(Y_T-\partial_x g(X_T, \mu^N_T)\big)^2\big]\\
    \text{s.t: }   &  X_{t_{n+1}} = X_{t_n} + \Delta t \times b(t_n, X_{t_n}, \hat \alpha_{t_n}, \mu^N_{t_n})+\sigma \Delta W_{t_n}, \quad X_0 \sim \mu_0 \\
    &Y_{t_{n+1}} = Y_{t_n} - \Delta t \times \partial_x H(t_n, X_{t_n}, \hat \alpha_{t_n}, \mu^N_{t_n})+Z_{t_n} \Delta W_{t_n},\ \forall t_n \in \mathcal{T} \backslash T,
    \end{aligned}
\end{equation}
where $\hat \alpha_{t_n} = \hat{\alpha}(t_n, X_{t_n}, \mu^N_{t_n}, Y_{t_n})$ and $\Delta W_{t_n}\sim \mathcal{N}(0,\Delta t)$ with $\Delta t = t_{n+1}-t_n$ where $t_n$ represents the discretized time steps in $\mathcal{T} = \{0, t_1, t_2, \dots, T\}$. In order to approximate the mean field, we  simulate by Monte Carlo a population of size $N$ and use the empirical distribution of states $\mu_{t_n}^N = \frac{1}{N}\sum_{i=1}^N\delta_{X_{t_n}^i}$ in the population. In this method, we can directly aim at minimizing the shooting error by replacing controls $(Y_0, \bZ)$ by parameterized functions $y_{0, \theta_1} : \RR \rightarrow \RR$ and $z_{\theta_2}: [0,T] \times \RR \rightarrow \RR$, in which, as in Section~\ref{sec:FBSDE-iterative-algo}, $Y_0$ is implemented as a function of $X_0$ and $\bZ$ is implemented as a function of $t$ and $X_t$. In order to approximate these functions we will use neural networks similar to the previous algorithm. As we discussed in the previous algorithm, we can minimize the cost in~\eqref{eq:shooting} by using SGD by sampling one population of size $N$. Precisely, the aim is to minimize over the neural network parameters $\theta = (\theta_1, \theta_2)$ the following (empirical) objective: 
\begin{equation}
\label{eq:averaged_shooting_XY}
\begin{aligned}
    \tilde{\mathbb{J}}^N(\theta) = \EE \left[ \dfrac{1}{N}\sum_{i=1}^N \big(Y^{i,\theta}_T-\partial_x g(X^{i, \theta}_T, \mu_T^{N, \theta})\big)^2 \right].
\end{aligned}
\end{equation}
The main difference of this algorithm from the previous one is that it simulates the trajectories of $\bX$ and $\bY$ simultaneously.
We can see that the objective to minimize in~\eqref{eq:averaged_shooting_XY} is different than the objective in~\eqref{eq:averaged_shooting_Y} since now $\bX$ and $\bmu$ are directly affected by the neural network parameters $\theta$. The details of the algorithm can be found in Algorithms~\ref{algo:simu-forward-XY} and~\ref{algo:SGD-MFG-Alg2}. 

This type of algorithm has been used e.g. by~\citet{fouque2020deep,carmona2022convergence,germain2022numerical,epidemics_SMFG,dayanikli2023machine} to solve various kinds of forward-backward mean field SDE systems.

    \begin{algorithm}[h]
    \caption{Simultaneous Monte Carlo simulation of an interacting batch $\bX$ and $\bY$ \label{algo:simu-forward-XY}}
    
    \textbf{Input:} number of particles $N$; time horizon $T$; time increment $\Delta t$; initial distribution $\mu_{0}$; control functions $ y_0, z$
    
    \textbf{Output:} Approximate sampled trajectories of $(\bX,\bY,\bZ)$.
    \begin{algorithmic}[1]
    \State{Let $n=0$, $t_0 = 0$; pick $X^i_0 \sim \mu_{0}$ i.i.d. and set $Y^i_0 = y_0(X^i_0), \forall i \in [N]$}
    \While{$n \times \Delta t = t_n < T$}
        \State{Set $ Z^i(t_n) = z({t_n}, X^i_{t_n}),\ i \in [N]$}\vskip1mm
        \State{Let $ \mu^N_{t_n} = \frac{1}{N} \sum_{i=1}^N \delta_{X^i_{t_n}}$ }\vskip1mm
        \State{Set $\hat \alpha^i_{t_n} = \hat{\alpha}({t_n},X^i_{t_n},{\mu}^N_{t_n}, Y^i_{t_n}),\ i \in [N]$}\vskip1mm
        \State{Sample $\Delta W^i_{t_n} \sim \mathcal{N}(0,\Delta t),\ i \in [N]$ }\vskip1mm
        \State{Let 
            $
                X^{i}_{t_{n+1}} = X^{i}_{t_n} + b(t_n, X^{i}_{t_n}, \hat\alpha^{i}_{t_n}, \mu^N_{t_n})\Delta t + \sigma \Delta W^{i}_{t_n}
            $, $ i \in [N]$}\vskip1mm
        \State{Let 
            $
                Y^{i}_{t_{n+1}} = Y^{i}_{t_n} - \partial_x H(t_n, X^{i}_{t_n}, \hat \alpha^{i}_{t_n}, \mu^N_{t_n})\Delta t + Z^{i}_{t_n} \Delta W^{i}_{t_n}
            $, $ i \in [N]$
        }\vskip1mm
        \State{Set $t_n = t_n+\Delta t$ and $n=n+1$}
    \EndWhile
    \State{Set $N_T = n,\ t_{N_T} = T$}\\
    \Return $(X^i_{t_n},Y^i_{t_n},Z^i_{t_n})_{n =0,\dots,N_T,\ i \in [N]}$ 
    \end{algorithmic}
    \end{algorithm}

\begin{algorithm}[h]
\caption{Stochastic Gradient Descent for solving~\eqref{eq:shooting} \label{algo:SGD-MFG-Alg2}}

\textbf{Input:} Initial parameter $\theta_0$; number of iterations ${K}$; sequence $(\beta_{{k}})_{{k} = 0, \dots, {K}-1}$ of learning rates; number of particles $N$; time horizon $T$; time increment $\Delta t$; initial distribution $\mu_0$

\textbf{Output:} Approximation of $\theta^*$ minimizing $\tilde{\mathbb{J}}^N$ defined by~\eqref{eq:averaged_shooting_XY}
\begin{algorithmic}[1]

\For{${k} = 0, 1, 2, \dots, {K}-1$}
    \State \parbox[t]{\dimexpr\textwidth-\leftmargin-\labelsep-\labelwidth}{Sample $S = (X^i_{t_{{n}}},Y^i_{t_{{n}}},Z^i_{t_{{n}}})_{{n} =0,\dots,{N_T},\ i \in  [N]}$ 
    using Algorithm~\ref{algo:simu-forward-XY} with control functions $(y_0, \boldsymbol z) = (y_{0,\theta_{{k},1}}, z_{\theta_{{k},2}})$ and parameters: number of particles $N$; time horizon $T$; time increment $\Delta t$; initial distribution $\mu_{0}$ \strut}
    \State{Compute the gradient $\nabla \tilde{\mathbb{J}}^N(\theta_{{k}})$ of $\tilde{\mathbb{J}}^N(\theta_{{k}})$}
    \State{Set $\theta_{{k}+1} =  \theta_{{k}} -\beta_{{k}} \nabla \tilde{\mathbb{J}}^N(\theta_{{k}})$ 
    }
\EndFor

\State \Return $\theta_{{K}}$
\end{algorithmic}
\end{algorithm}

\subsubsection{Numerical illustration: Epidemic model with FBSDE deep learning}
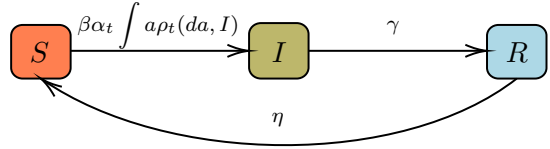
\begin{figure}[H]
\begin{center}

\tikzset{every picture/.style={line width=0.75pt}}
\begin{tikzpicture}[x=0.75pt,y=0.7pt,yscale=-1,xscale=1]

%Rounded Rect [id:dp5932690496195532] susceptible
\draw  [fill={rgb, 255:red, 255; green, 127; blue, 80 }  ,fill opacity=1 ] (100,106.8) .. controls (100,103.6) and (102.6,101) .. (105.8,101) -- (124.2,101) .. controls (127.4,101) and (130,103.6) .. (130,106.8) -- (130,124.2) .. controls (130,127.4) and (127.4,130) .. (124.2,130) -- (105.8,130) .. controls (102.6,130) and (100,127.4) .. (100,124.2) -- cycle ;
%Rounded Rect [id:dp5633636707902099] infected
\draw  [fill={rgb, 255:red, 189; green, 183; blue, 107 }  ,fill opacity=1 ] (220,106) .. controls (220,102.69) and (222.69,100) .. (226,100) -- (244,100) .. controls (247.31,100) and (250,102.69) .. (250,106) -- (250,124) .. controls (250,127.31) and (247.31,130) .. (244,130) -- (226,130) .. controls (222.69,130) and (220,127.31) .. (220,124) -- cycle ;
%Rounded Rect [id:dp5169215905527198] recovered
\draw  [fill={rgb, 255:red, 173; green, 216; blue, 230 }  ,fill opacity=1 ] (340,106.8) .. controls (340,103.6) and (342.6,101) .. (345.8,101) -- (364.2,101) .. controls (367.4,101) and (370,103.6) .. (370,106.8) -- (370,124.2) .. controls (370,127.4) and (367.4,130) .. (364.2,130) -- (345.8,130) .. controls (342.6,130) and (340,127.4) .. (340,124.2) -- cycle ;
%Straight Lines [id:da7643207816793263] 
\draw    (130,115) -- (218,115) ;
\draw [shift={(220,115)}, rotate = 180] [color={rgb, 255:red, 0; green, 0; blue, 0 }  ][line width=0.75]    (10.93,-3.29) .. controls (6.95,-1.4) and (3.31,-0.3) .. (0,0) .. controls (3.31,0.3) and (6.95,1.4) .. (10.93,3.29)   ;
%Straight Lines [id:da2928359833861143] 
\draw    (250,115) -- (320.5,115) -- (338,115) ;
\draw [shift={(340,115)}, rotate = 180] [color={rgb, 255:red, 0; green, 0; blue, 0 }  ][line width=0.75]    (10.93,-3.29) .. controls (6.95,-1.4) and (3.31,-0.3) .. (0,0) .. controls (3.31,0.3) and (6.95,1.4) .. (10.93,3.29)   ;
%Curve Lines [id:da5849227676643367] 
\draw    (355,130) .. controls (296.79,178.76) and (178.69,177.02) .. (115.94,130.7) ;
\draw [shift={(115,130)}, rotate = 396.94] [color={rgb, 255:red, 0; green, 0; blue, 0 }  ][line width=0.75]    (10.93,-3.29) .. controls (6.95,-1.4) and (3.31,-0.3) .. (0,0) .. controls (3.31,0.3) and (6.95,1.4) .. (10.93,3.29)   ;

% Text Node
\draw (108,109) node [anchor=north west][inner sep=0.75pt]   [align=left] {$\displaystyle S$};
% Text Node
\draw (230,109) node [anchor=north west][inner sep=0.75pt]   [align=left] {$\displaystyle I$};
% Text Node
\draw (348,109) node [anchor=north west][inner sep=0.75pt]   [align=left] {$\displaystyle R$};
% Text Node
\draw (132,87) node [anchor=north west][inner sep=0.75pt]  [font=\scriptsize] [align=left] {$\displaystyle \beta \alpha _{t}\int a\rho _{t}( da,I)$};
% Text Node
\draw (288,97) node [anchor=north west][inner sep=0.75pt]  [font=\scriptsize] [align=left] {$\displaystyle \gamma $};
% Text Node
\draw (230,147) node [anchor=north west][inner sep=0.75pt]  [font=\scriptsize] [align=left] {$\displaystyle \eta$};

\end{tikzpicture}

\end{center}
\caption{SIR model with mean-field interactions.}
\label{fig:SIR-diagram-intro}
\end{figure}

\begin{figure}[h]
\centering
\begin{subfigure}{.32\textwidth}
    \includegraphics[width=4.5cm,height=3cm]{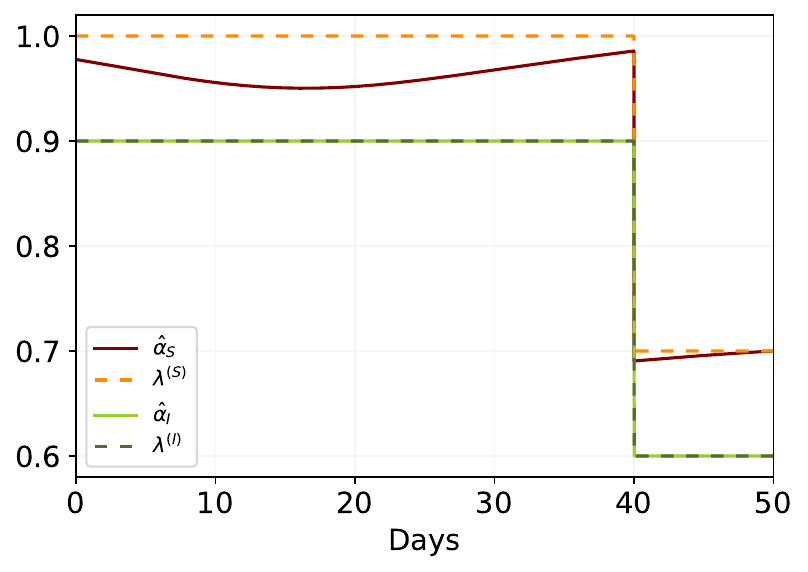}
\end{subfigure}
\begin{subfigure}{.32\textwidth}
\includegraphics[width=4.5cm,height=3cm]{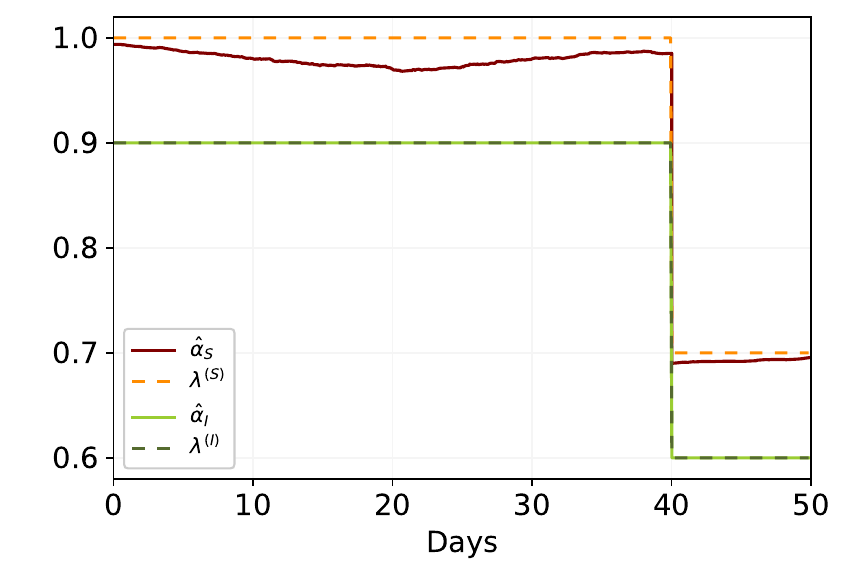}
\end{subfigure}
\begin{subfigure}{.32\textwidth}
\includegraphics[width=4.9cm,height=3.3cm]{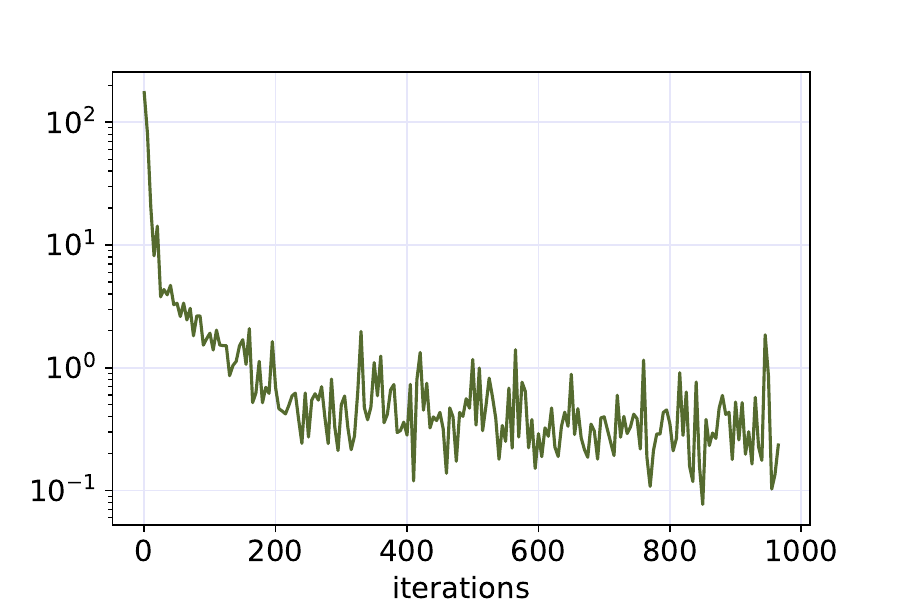}
\end{subfigure}
\caption{Evolution of the controls: explicit solution (left), deep learning solution (middle), convergence of the loss value in the deep learning approach (right). }
\label{fig:LL-NN}
\end{figure}

Before concluding this section, we introduce an example of mitigating epidemics from~\citet*{epidemics_SMFG} that solves FBSDEs by using deep learning. Different than the setup introduced here, this model is a \textit{finite} state model. However, for showing the general applicability of the deep learning approach, we introduce the example and the results. The representative agent can be in 3 different states: Susceptible ($S$), Infected ($I$), and Recovered ($R$). The representative agent's state transitions from $S$ to $I$ depending on her own socialization level $\alpha_t$, exogenous infection rate $\beta>0$, and the average socialization level of the infected people $\int a\rho _{t}(da,I)$, where $\rho$ is the joint control and state distribution. The last term introduces the mean field interactions into the model. Here, the socialization level $\alpha_t$ is the action chosen by the representative agent at time $t$. The interactions are thus through the distribution of actions. The representative agent's state transitions from $I$ to $R$ and $R$ to $S$ with constant rates $\gamma$ and $\kappa$ respectively. The jump rates can be seen in Figure~\ref{fig:SIR-diagram-intro}. We consider that the objective of the representative agent is to minimize the following cost while choosing their socialization level $\balpha=(\alpha_t)_{t\in[0,T]}$:
\begin{equation}
\label{eq:cost-agents-SIR1-intro}
    \EE\Big[\int_0^T 
    \tfrac{c_\lambda \left|\lambda^{(S)}_t - \alpha_t\right|^2}{2}\mathbbm{1}_{S}(x) 
    + \Big(\tfrac{\left|\lambda^{(I)}_t - \alpha_t\right|^2}{2} + c_I\Big)\mathbbm{1}_{I}(x) 
    + \tfrac{\left|\lambda^{(R)}_t-\alpha_t\right|^2}{2}\mathbbm{1}_{R}(x) dt\Big],
\end{equation}
where $c_\lambda$ and $c_I$ are positive exogenous constants. In the cost~\eqref{eq:cost-agents-SIR1-intro}, the terms $\lambda\in[0,1]$ represent the social distancing policies (such as mask wearing, quarantining) set by the government exogenously and individuals do not want to deviate from these policies significantly. Moreover, infected people have a treatment cost denoted by $c_I$. We borrow the results from~\citet{epidemics_SMFG} and show them in Figure~\ref{fig:LL-NN}. On the left plot and middle plot of Figure~\ref{fig:LL-NN}, we can see the social distancing policies set by the government with the dashed lines (orange for susceptible individuals and green for infected individuals) and the mean field Nash equilibrium socialization levels of the susceptible (maroon) and infected (dark green) individuals. On the $40^{\rm th}$ day, the government is setting more restrictive social distancing policies (which is implemented in the model by setting lower $\lambda$ values). At the equilibrium, Infected individuals follow the policies closely, but susceptible individuals do not follow the policies set by the government closely instead they restrict their socialization level further to protect their health. In this figure the left plot corresponds to the \textit{explicit} solution and the middle plot corresponds to the \textit{deep learning} approach implemented with the algorithm of Section~\ref{sec:fbsde-algo2-simultaneous} above. We can see the deep learning method obtains a behavior which is very similar to the explicit solution.

\subsection{Deep learning methods for solving PDEs}
\label{subsec:deeplearning_fbpde}
MFG equilibrium can be also characterized by using forward backward partial differential equations (FBPDEs). In this formulation, the backward equation is the Hamilton-Jacobi-Bellman (HJB) equation and it characterizes the value function defined in Equation~\eqref{eq:value}. The forward equation is the Kolmogorov-Fokker-Planck (KFP) equation and it characterizes the flow of the mean field. Then, the equilibrium control is given as the minimizer of the Hamiltonian, i.e., 
$ \hat{\alpha}(t,x) =\argmin_{\alpha} H(t, x, \alpha, m(t), \partial_x u(t, x))$, where Hamiltonian is defined as follows
\begin{equation*}
    H(t, x, \alpha, m, p) = b(t, x, \alpha, m) p + f(t, x, \alpha, m), 
\end{equation*}
and $(\bm, \bu)$ solves the following FBPDE:
\begin{equation}
\label{eq:fbpde_general}
    \begin{aligned}
        \partial_t{m(t,x)} &= \frac{1}{2}\sigma^2 \partial^2_{xx} m(t,x) - \partial_x (b(t, x, \hat{\alpha}(t,x), m_t)m(t,x)), && m(0,x)=m_0(x) \\
        -\partial_t u(t,x) &= \frac{1}{2}\sigma^2 \partial^2_{xx} u(t,x) + H(t, x, \hat{\alpha}(t,x), m_t, \partial_x u(t,x)), && u(T,x)=g(x, m_T).
    \end{aligned}
\end{equation}
Using the convention from the literature, we replaced the notation $\mu_t$ with $m_t$.

\subsubsection{Deep Galerkin method for mean field FBPDEs}

Similarly to the FBSDE approach discussed above, in the FBPDE approach, the backward and the forward equations are coupled and cannot be solved individually. In order to solve them, we use an extension to the Deep Galerkin method (DGM) introduced by~\citet{sirignano2018dgm}. DGM is a method that is introduced to solve general PDEs which shares similarities with physics informed neural networks~\citet{karniadakis2021physics}. In this section we explain how it can be implemented to solve the FBPDEs that characterize a solution of the MFGs. Interested reader can refer to~\citet{al2018solving,cao2020connecting,ruthotto2020machine,carmona2021convergence} for further implementations of DGM to solve FBPDEs.

For the implementation of this method, we will replace the state space with a compact subset $\mathcal{R}$ of $\RR$. We use neural networks to approximate the density function $m(t,x)$ and the value function $u(t,x)$ by rewriting the FBPDE as an optimization problem in which we aim to minimize a loss function. This loss is constituted by the PDE \textit{residuals} and penalty terms for the boundary conditions. It can be written as follows:
\begin{equation}
\label{eq:fbpde_cost}
    \mathbb{L}(\bm,\bu) = L^{f} (\bm,\bu) + L^{b} (\bm,\bu),
\end{equation}
where the $L^{f}(m,u)$ denotes the costs related to the forward component (i.e., KFP) and $L^{b}(m,u)$ denotes the costs related to the backward component (i.e., HJB). These cost components are defined as follows:
\begin{align}
    L^{f} (\bm,\bu) = &c_f \lVert \partial_t{m} - \frac{1}{2}\sigma^2 \partial^2_{xx} m + \partial_x (b(t, x, \hat{\alpha}, m(t))m)\rVert_{L^2([0,T] \times \mathcal{R})}+  c_i  \lVert m(0)-m_0\rVert_{L^2(\mathcal{R})},\\
    L^{b} (\bm,\bu) = &c_b \lVert \partial_t u +\frac{1}{2}\sigma^2 \partial^2_{xx} u + H(t, x, \hat{\alpha}, m(t), \partial_x u)\rVert_{L^2([0,T] \times \mathcal{R})}+c_t  \lVert u(T)-g(\cdot, m_T)\rVert_{L^2(\mathcal{R})},
\end{align}
where $\lVert \phi \rVert_{L^2(A)}$ defines the $L_2$-norm of function $\phi$ defined on space $A$, and $c_f, c_i, c_b, c_t$ are positive constants which are chosen to balance the importance of the corresponding terms in the total cost. If a density and value couple $(\hat{\bm}, \hat{\bu})$ solves the FBPDE~\eqref{eq:fbpde_general}, then we can see that the corresponding cost $\mathbb{L}(\hat{\bm},\hat{\bu})$ should be zero. The deep learning approach aims to minimize this cost by replacing functions $(\bm, \bu)$ by parameterized functions $m_{\theta_1}: [0, T] \times \mathcal{R} \to [0,1]$ and $u_{\theta_2}: [0, T] \times \mathcal{R} \to \RR$ where $\theta=(\theta_1, \theta_2)$ are the neural network parameters. As we introduced in Section~\ref{subsec:deeplearning_fbsde}, we can use SGD to minimize the cost~\eqref{eq:fbpde_cost} over the neural network parameters $\theta$. Here, this means the following. To apply mini-batch gradient descent, at each iteration we sample $\bS = (S, S_i, S_t)$ where $S$ is a finite set of points in $[0,T]\times\mathcal{R}$, and $S_i$, $S_t$ are two finite sets of points in $\mathcal{R}$. Then, the (empirical) cost function is computed as:
\begin{equation}
\label{eq:fbpde_cost_sample}
    \mathbb{L}^{\bS}(\theta) = L^{\bS,f} (\theta) + L^{\bS,b} (\theta),
\end{equation}
where $L^{\bS,f} = c_f L^{\bS,f,\mathrm{res}} + c_i L^{\bS,f,0}$ and $L^{\bS,b} = c_b  L^{\bS,b,\mathrm{res}} + c_t L^{\bS,b,T}$, with:
\begin{align}
    L^{\bS,f,\mathrm{res}} (\theta) &= \dfrac{1}{|S|} \sum_{(t,x)\in S} \big|\partial_t{m_{\theta_1} (t,x)} - \frac{1}{2}\sigma^2 \partial^2_{xx} m_{\theta_1}(t,x) + \partial_x \big(\hat{b}_\theta(t, x, m_{\theta_1}(t))m_{\theta_1}(t,x)\big)\big|^2,
    \\
     L^{\bS,f,0} (\theta) &=  \dfrac{1}{|S_i|} \sum_{x\in S_i}\big| m_{\theta_1}(0,x)-m_0(x)\big|^2,
     \\
    L^{\bS,b,\mathrm{res}} (\theta) = &\dfrac{1}{|S|} \sum_{(t,x)\in S} \big| \partial_t u_{\theta_2}(t,x) +\frac{1}{2}\sigma^2 \partial^2_{xx} u_{\theta_2}(t,x) + \hat{H}_\theta(t, x, m_{\theta_1}(t), \partial_x u_{\theta_2}(t,x))\big|^2,
    \\
     L^{\bS,b,T} &=  \dfrac{1}{|S_t|} \sum_{x\in S_t}\big| u_{\theta_2}(T,x)-g(x, m_{\theta_1}(T))\big|^2,
\end{align}
where we used the shorthand notations $\hat{b}_\theta(t, x, m_{\theta_1}(t)) = b(t, x, \hat{\alpha}_{\theta}(t,x), m_{\theta_1}(t))m_{\theta_1}(t,x)$ and\\ $\hat{H}_\theta(t, x, m_{\theta_1}(t), \partial_x u_{\theta_2}(t,x)) = H(t, x, \hat{\alpha}_{\theta}(t,x), m_{\theta_1}(t), \partial_x u_{\theta_2}(t,x))$ with the control function $\hat{\alpha}_{\theta}(t,x) = \hat{\alpha} (t,x,m_{\theta_1}(t),\partial_x u_{\theta_2}(t,x))$, which depends on both components of $\theta = (\theta_1,\theta_2)$.
The details of the algorithm can be found in Algorithm~\ref{algo:DGM-MFG}.

\begin{algorithm}[H]
\caption{Deep Galerkin Method for solving FBPDE system~\eqref{eq:fbpde_general} \label{algo:DGM-MFG}}

\textbf{Input:} Initial parameter $\theta_0$; number of iterations ${K}$; sequence $(\beta_{{k}})_{{k} = 0, \dots, {K}-1}$ of learning rates.

\textbf{Output:} Approximation of $\theta^*$ minimizing $\mathbb{L}^{\bS}$ defined in~\eqref{eq:fbpde_cost_sample}
\begin{algorithmic}[1]

\For{${k} = 0, 1, 2, \dots, {K}-1$}
    \State \parbox[t]{\dimexpr\textwidth-\leftmargin-\labelsep-\labelwidth}{Sample $\bS = (S, S_i, S_t)$\strut}
    \State{Compute the gradient $\nabla \mathbb{L}^{\bS}(\theta_{{k}})$ of $\mathbb{L}^{\bS}(\theta_{k})$}
    \State{Set $\theta_{{k}+1} =  \theta_{{k}} -\beta_{{k}}  \nabla \mathbb{L}^{\bS}(\theta_{{k}})$ }
\EndFor

\State \Return $\theta_{{K}}$
\end{algorithmic}
\end{algorithm}

\subsubsection{Numerical illustration: Portfolio liquidation with PDE deep learning} 
\label{sec:numerics-trading}
We consider the model proposed by~\citet{MR3805247} and discussed in Section~\ref{sec:example trading}. Following~\citet{MR3805247},  the Nash equilibrium control is $\alpha^*_t(q) = \frac{\partial_q v(t,q)}{2\kappa}$, where $(v,m)$ solve the following PDE system:
\begin{equation}
\label{eq:trading-PDE-system}
\left\{
\begin{aligned}
    &\quad -\gamma\bar{\mu} q = \partial_t v - \phi q^2 + \frac{|\partial_q v(t,q)|^2}{4\kappa}
    \\
    &\quad \partial_t m + \partial_q\left(m \frac{\partial_q v(t,q)}{2\kappa}\right) = 0
    \\
    &\quad \bar{\mu}_t = \int \frac{\partial_q v(t,q)}{2\kappa} m(t, dq)
    \\
    &\quad m(0,\cdot) = m_0, v(T,q) = - A q^2. 
\end{aligned}
\right.
\end{equation}
As already mentioned before, the interactions are non-local and involve the mean of the controls and not just the mean of the states.

\begin{figure}[hb]
\centering
\begin{subfigure}{.45\textwidth}
  \centering
  \includegraphics[width=\linewidth]{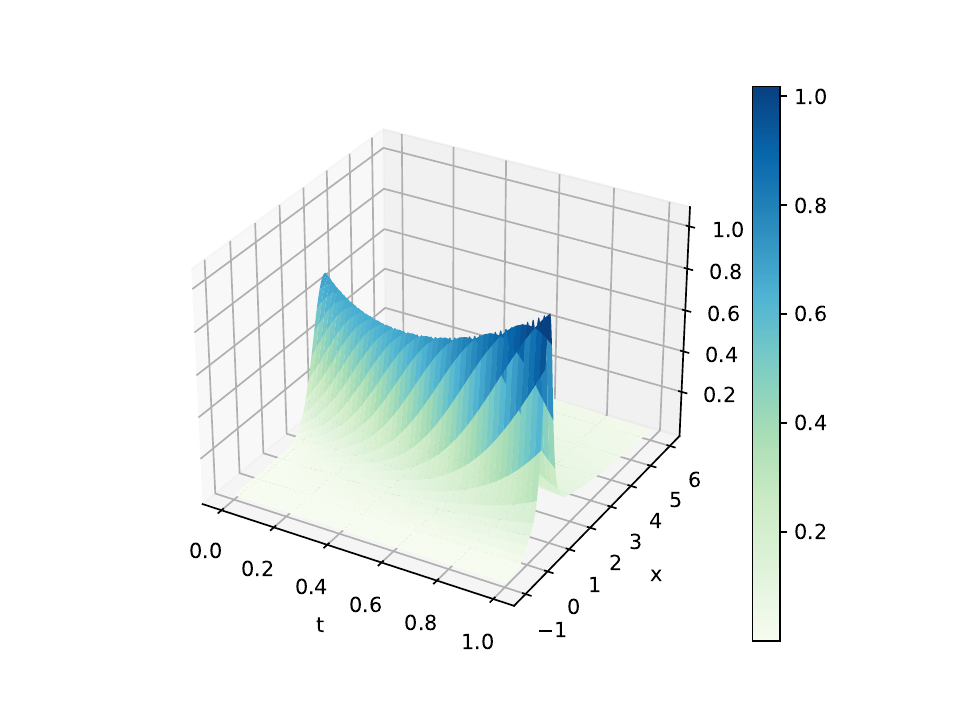}
\end{subfigure}
\begin{subfigure}{.45\textwidth}
  \centering
  \includegraphics[width=\linewidth]{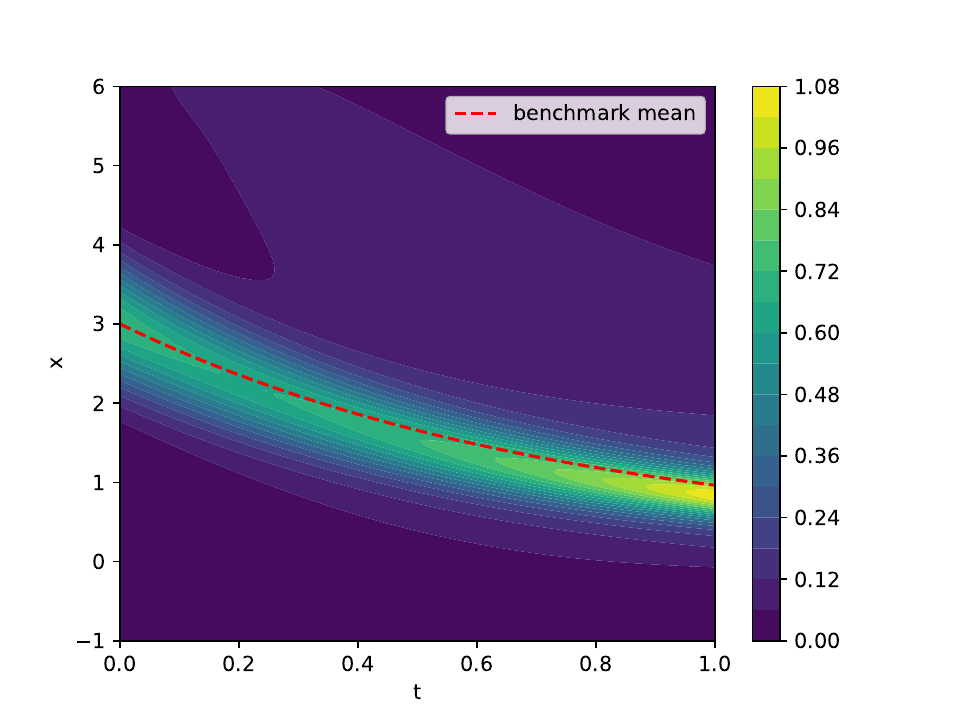}
\end{subfigure}
\caption{Case 1 of trading MFG example solved by DGM. Evolution of the distribution $m$: surface (left) and contour (right). The dashed red line corresponds to the mean obtained by the semi-explicit formula. }
\label{fig:ex-mfg-trading-sell-distrib}
\end{figure}

\begin{figure}[ht]
\centering
\begin{subfigure}{.3\textwidth}
  \centering
  \includegraphics[width=\linewidth]{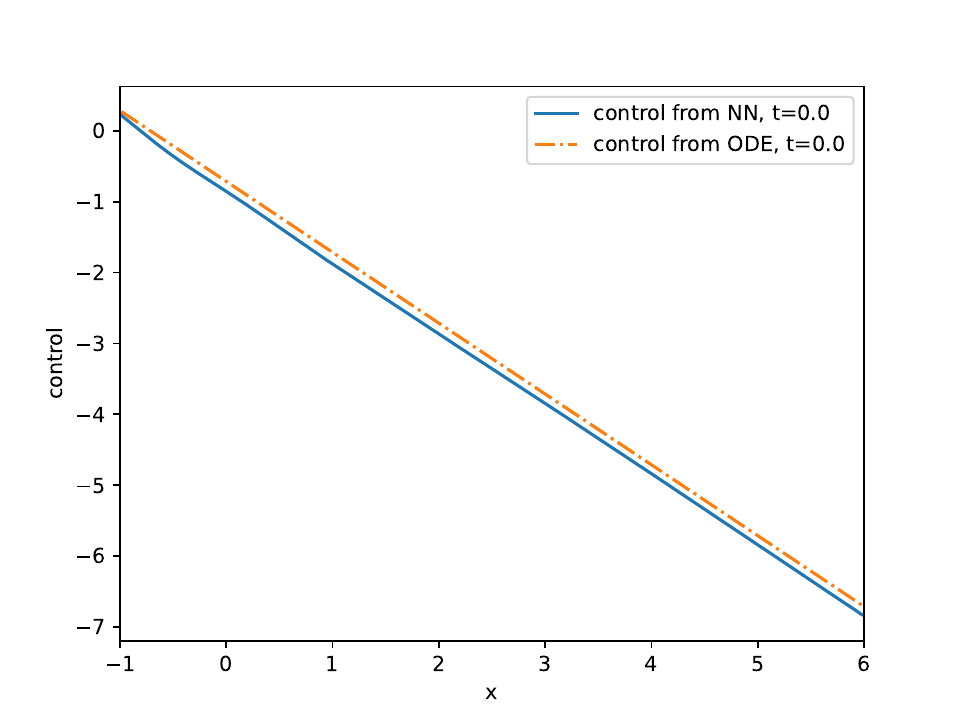}
\end{subfigure}
\begin{subfigure}{.3\textwidth}
  \centering
  \includegraphics[width=\linewidth]{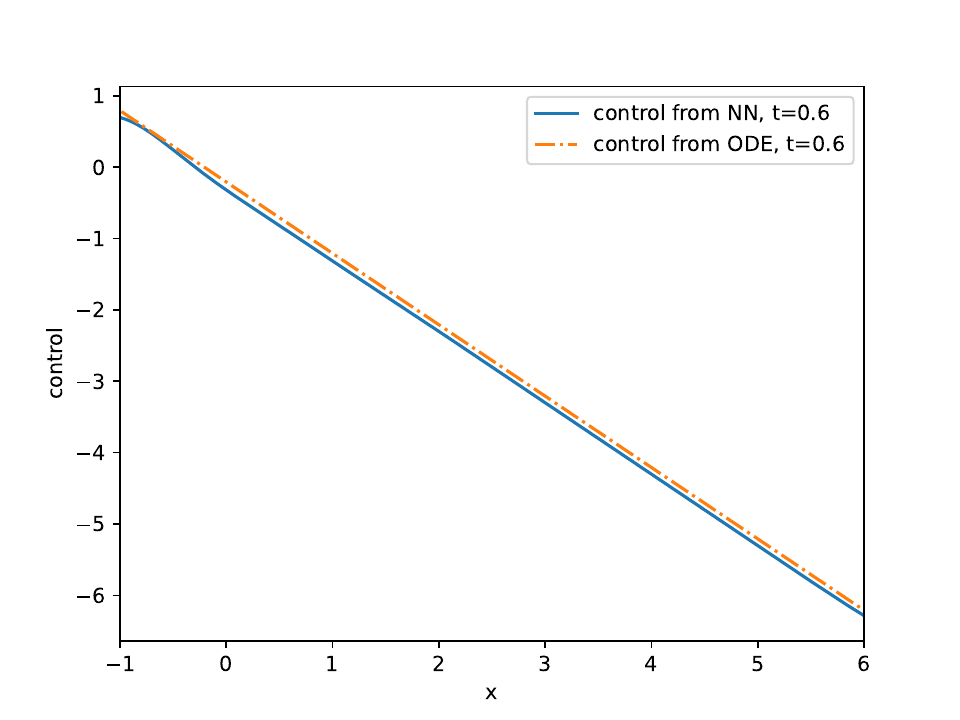}
\end{subfigure}
\begin{subfigure}{.3\textwidth}
  \centering
  \includegraphics[width=\linewidth]{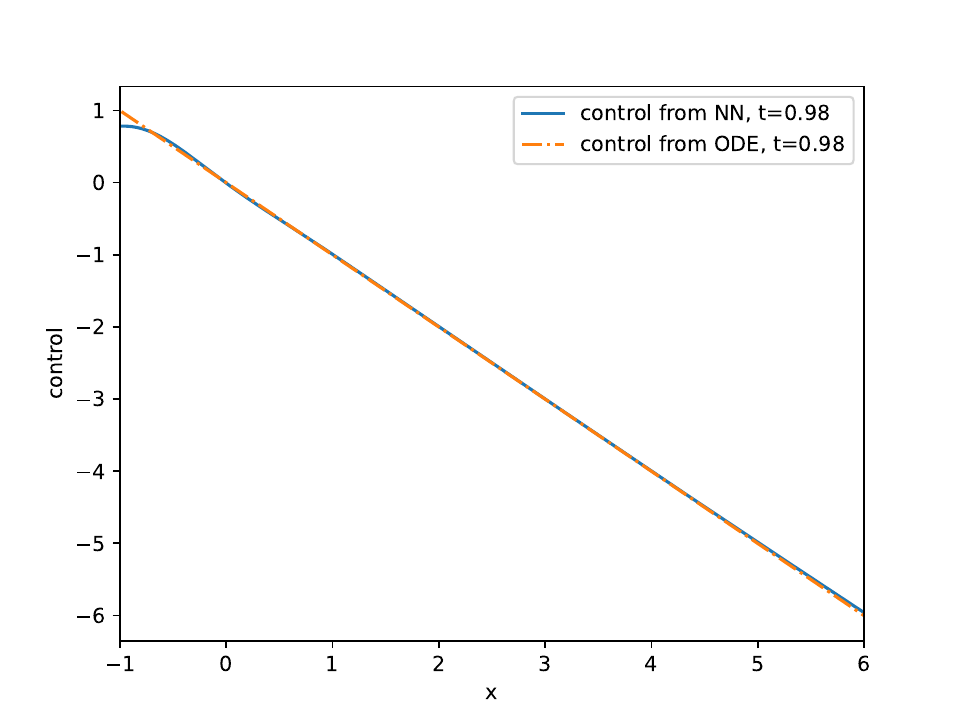}
\end{subfigure}
\caption{Case 1 of trading MFG example solved by DGM. Optimal control $\alpha^*$ (dashed line) and learnt control (full line)  at three different time steps.}
\label{fig:ex-mfg-trading-sell-ctrl}
\end{figure}

\begin{figure}[ht]
\centering
\begin{subfigure}{.45\textwidth}
  \centering
  \includegraphics[width=\linewidth]{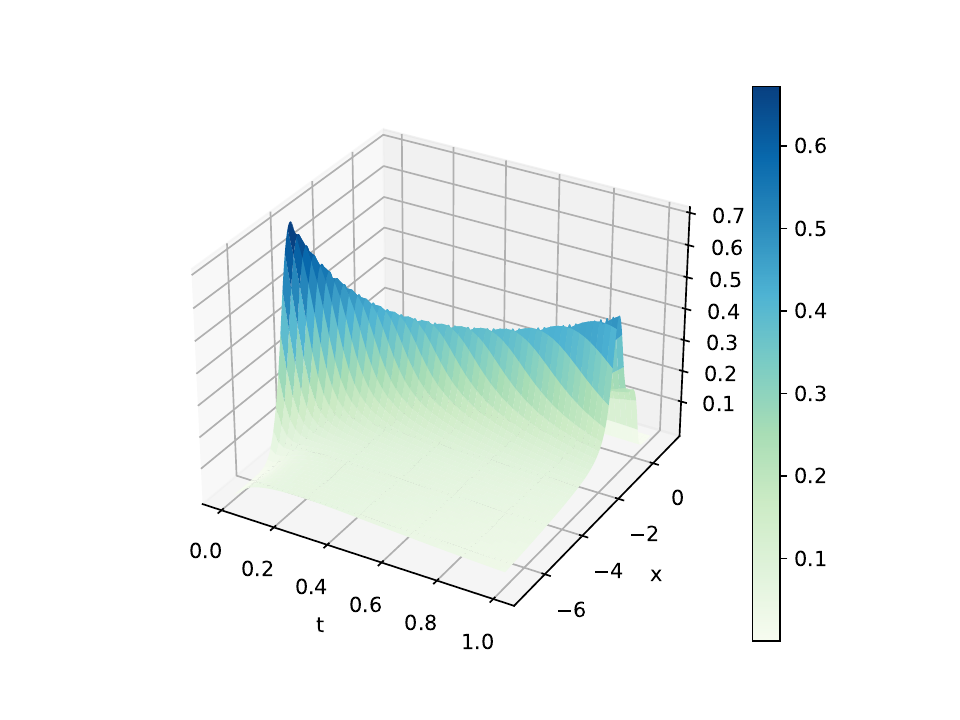}
\end{subfigure}
\begin{subfigure}{.45\textwidth}
  \centering
  \includegraphics[width=\linewidth]{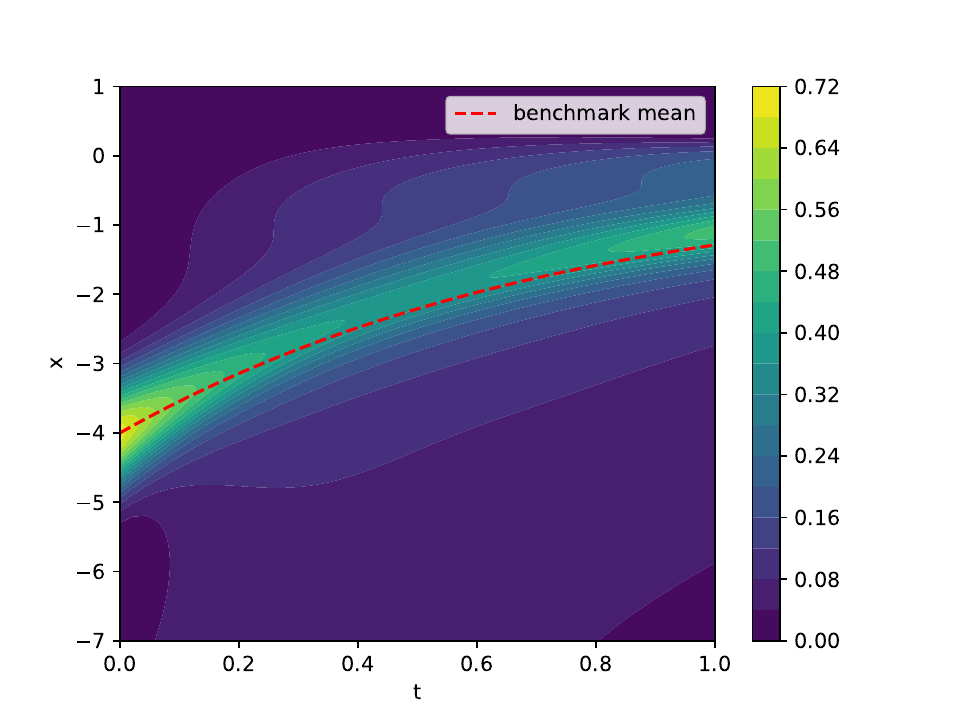}
\end{subfigure}
\caption{Case 2 of trading MFG example solved by DGM. Evolution of the distribution $m$: surface (left) and contour (right). The dashed red line corresponds to the mean obtained by the semi-explicit formula. }
\label{fig:ex-mfg-trading-buy-distrib}
\end{figure}

\begin{figure}[ht]
\centering
\begin{subfigure}{.3\textwidth}
  \centering
  \includegraphics[width=\linewidth]{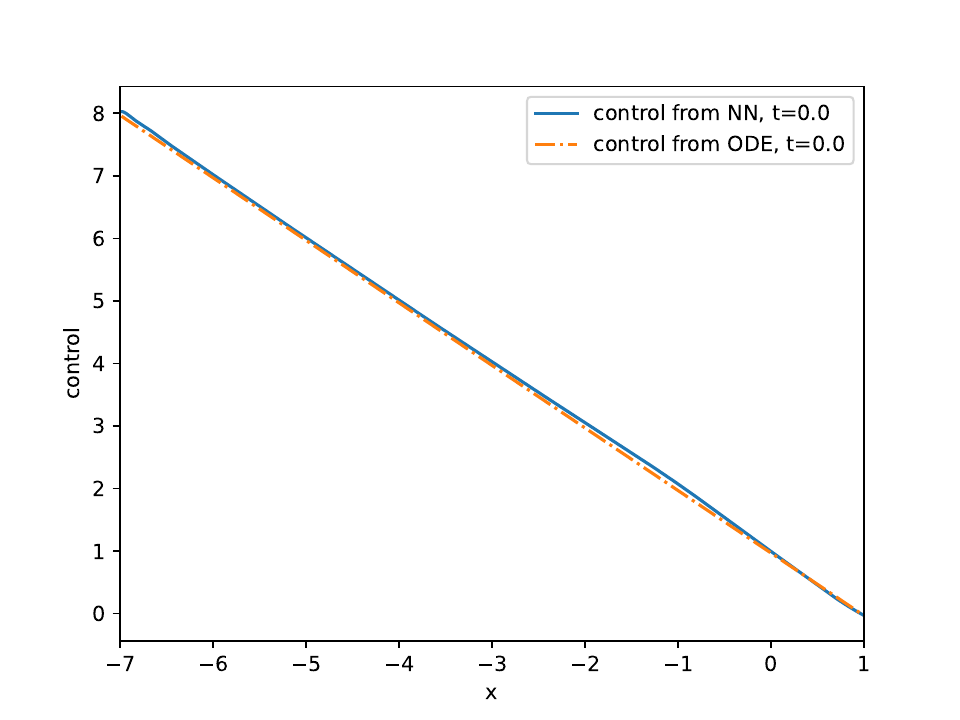}
\end{subfigure}
\begin{subfigure}{.3\textwidth}
  \centering
  \includegraphics[width=\linewidth]{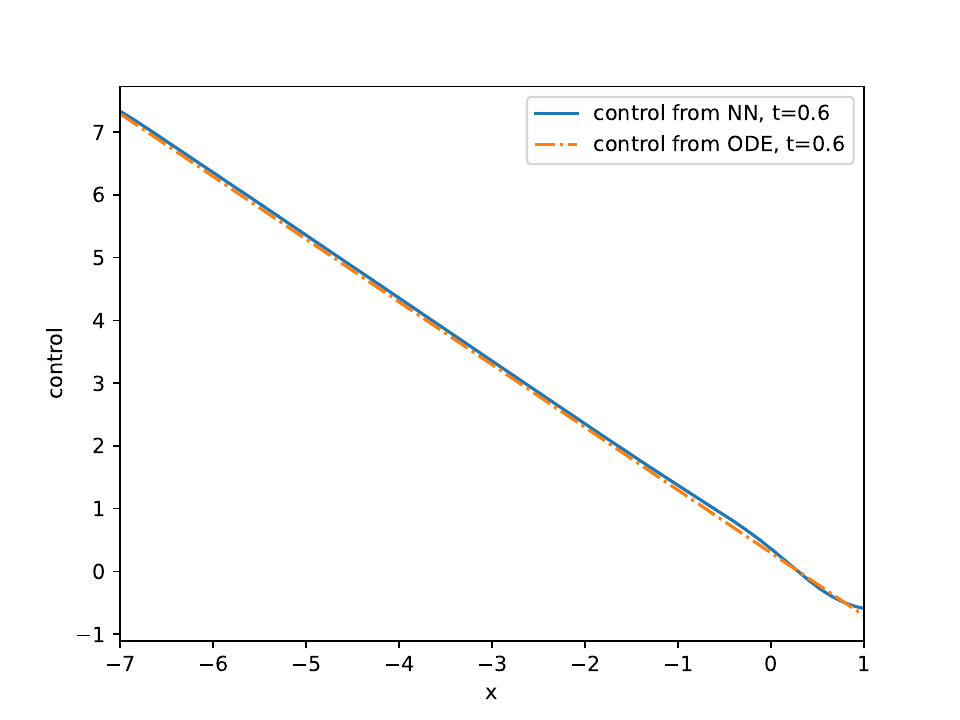}
\end{subfigure}
\begin{subfigure}{.3\textwidth}
  \centering
  \includegraphics[width=\linewidth]{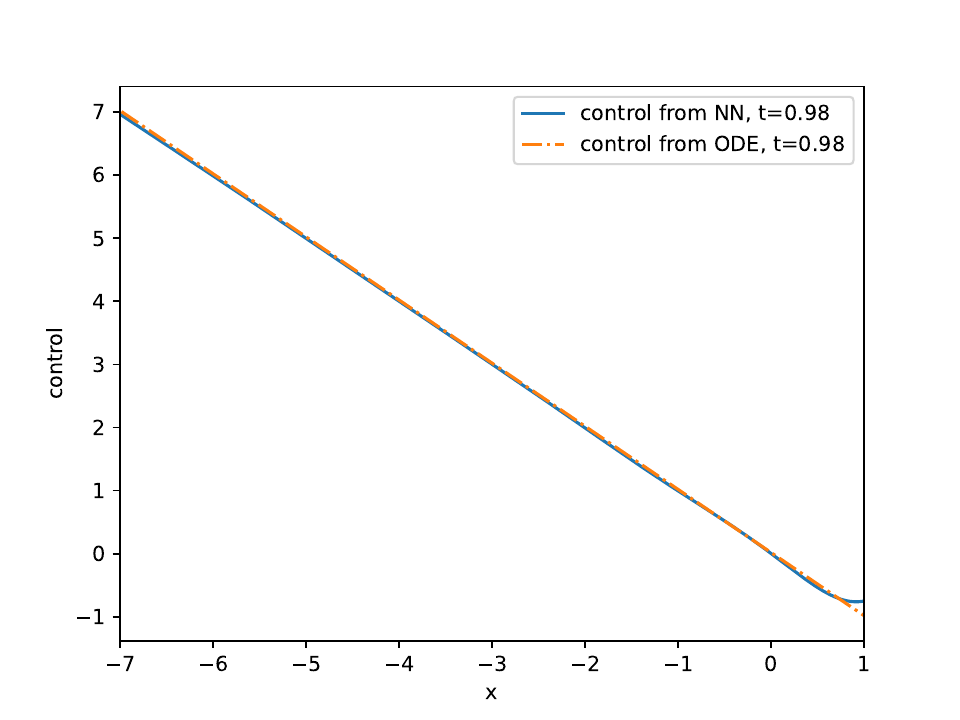}
\end{subfigure}
\caption{Case 2 of trading MFG example solved by DGM. Optimal control $\alpha^*$ (dashed line) and learnt control (full line)  at three different time steps.}
\label{fig:ex-mfg-trading-buy-ctrl}
\end{figure}

We solve this PDE system using the DGM.
We use the architecture proposed in the DGM article with 3 layers and a width of 40.
In the following results, we used the parameters: $T=1$, $A = 1$, $\phi = 1$, $\kappa = 1$, $\gamma = 1$, and a Gaussian initial distribution. We consider two cases: mean $3$ and variance $0.5$, and mean $-4$ and variance $0.3$. In both cases, the traders want to bring their inventory close to zero (but perhaps not exactly equal due to other terms in the cost and due to the price impact). In the first case, it means selling stocks, while in the second case, it means buying stocks. 
In the first case, Figures~\ref{fig:ex-mfg-trading-sell-distrib} and~\ref{fig:ex-mfg-trading-sell-ctrl} show respectively the evolution of the distribution and of the control, and the comparison with the solution coming from a benchmark obtained using ODEs.\footnote{In this model, the equilibrium can be reduced to solving a system of ODEs; see~\citet{MR3805247}. } In the second case, Figures~\ref{fig:ex-mfg-trading-buy-distrib} and~\ref{fig:ex-mfg-trading-buy-ctrl} show respectively the distribution and the controls.

\section{Conclusion and perspectives}
\label{sec:conclusion}

In the sections above, we have presented the basic ideas behind learning methods for large-population games, with a particular focus on mean field games. We discussed both discrete time and continuous time models.

From here, several directions remain to be investigated in future work.
For example, in the direction of real-world applications, it would be interesting to solve games -- and particularly MFGs --  with very complex and large environments. 
It would also be important to study how realistic the MFG solution is by comparing it with observed data and other solution methods such as agent-based models. Typical examples in OR could be in traffic routing, epidemics, electricity consumption or portfolio management. 
Another aspect that would make MFG models more realistic is to blend cooperative and non-cooperative behaviors, since in reality people do not behave purely in one of the two extreme cases. 
In the direction of the theoretical foundations, it would be interesting to establish the convergence of the algorithms presented in the previous sections under more flexible assumptions than what is known thus far. Last but not least, although deep neural networks yield strong empirical successes, conditions under which they break the curse of dimensionality remain little understood in the context of large population games.

\section*{Acknowledgements}

M.L. would like to thank Theophile Cabannes, Sertan Girgin, Julien P\'erolat and Kai Shao for their help with OpenSpiel for the examples displayed in this paper. The authors would like to thank the anonymous reviewers for their valuable comments.

% %%%%%%%%%%%%%%%%%%%%%%%%%%%%%%%%%%%%%%%%%%%%%%%%%%%%%
% \begin{APPENDIX}{}%{With a Title}
% %%%%%%%%%%%%%%%%%%%%%%%%%%%%%%%%%%%%%%%%%%%%%%%%%%%%%

% \end{APPENDIX}

\end{document}